\documentclass[reqno]{amsart}
\usepackage{amsaddr}
\usepackage{slashed}
\usepackage{enumerate}
\usepackage{mhequ}
\usepackage[margin=1.2in]{geometry}
\usepackage{amssymb,url,color, booktabs,nccmath}
\usepackage{mathrsfs}
\usepackage{enumitem}
\usepackage{graphicx}
\usepackage{tikz}
\usepackage{mhequ}
\usepackage{mhsymb}
\usepackage{longtable}
\usetikzlibrary{shapes,snakes}
\usetikzlibrary{calc}
\usetikzlibrary{decorations.shapes}
\usepackage{color}
\usepackage[colorlinks=true]{hyperref}
\hypersetup{
    linkcolor=blue,          
    citecolor=red,        
    filecolor=blue,      
    urlcolor=cyan
}

\usepackage{color}
\usepackage[colorlinks=true]{hyperref}
\hypersetup{
    linkcolor=blue,          
    citecolor=red,        
    filecolor=blue,      
    urlcolor=cyan
}

\usetikzlibrary{shapes.misc}
\usetikzlibrary{shapes.symbols}
\usetikzlibrary{snakes}
\usetikzlibrary{decorations}
\usetikzlibrary{decorations.markings}

\def\drawx{\draw[-,solid] (-3pt,-3pt) -- (3pt,3pt);\draw[-,solid] (-3pt,3pt) -- (3pt,-3pt);}
\tikzset{
	root/.style={circle, draw=black, fill=white, inner sep=0pt, minimum size=0.7mm},
	dot/.style={circle,fill=black,draw=black, solid,inner sep=0pt,minimum size=0.5mm},
	square/.style={rectangle,fill=black,draw=black, solid,inner sep=0pt,minimum size=1mm},
	empty/.style={circle,fill=white,draw=white, solid,inner sep=0pt,minimum size=0.5mm},
	var/.style={circle,fill=black!10,draw=black,inner sep=0pt, minimum size=
	2mm},
	symb/.style={circle,fill=symbols,draw=symbols, solid,inner sep=0pt,minimum size=0.5mm},
	yy/.style={circle,fill=gray!20,draw=black,inner sep=0pt,minimum size=0.8mm},
	>=stealth,
	dotred/.style={circle,fill=black!50,inner sep=0pt, minimum size=2mm},
	generic/.style={semithick,shorten >=1pt,shorten <=1pt},
	dist/.style={ultra thick,draw=testcolor,shorten >=1pt,shorten <=1pt},
	testfcn/.style={ultra thick,testcolor,shorten >=1pt,shorten <=1pt,<-},
	testfcnx/.style={ultra thick,testcolor,shorten >=1pt,shorten <=1pt,<-,
		postaction={decorate,decoration={markings,mark=at position 0.6 with {\drawx}}}},
	kprime/.style={semithick,shorten >=1pt,shorten <=1pt,densely dashed,->},
	kprimex/.style={semithick,shorten >=1pt,shorten <=1pt,densely dashed,->,
		postaction={decorate,decoration={markings,mark=at position 0.4 with {\drawx}}}},
	kernel/.style={semithick,shorten >=1pt,shorten <=1pt,->},
	multx/.style={shorten >=1pt,shorten <=1pt,
		postaction={decorate,decoration={markings,mark=at position 0.5 with {\drawx}}}},
	kernelx/.style={semithick,shorten >=1pt,shorten <=1pt,->,
		postaction={decorate,decoration={markings,mark=at position 0.4 with {\drawx}}}},
	kernel1/.style={->,semithick,shorten >=1pt,shorten <=1pt,postaction={decorate,decoration={markings,mark=at position 0.45 with {\draw[-] (0,-0.1) -- (0,0.1);}}}},
	kernel2/.style={->,semithick,shorten >=1pt,shorten <=1pt,postaction={decorate,decoration={markings,mark=at position 0.45 with {\draw[-] (0.05,-0.1) -- (0.05,0.1);\draw[-] (-0.05,-0.1) -- (-0.05,0.1);}}}},
	kernelBig/.style={semithick,shorten >=1pt,shorten <=1pt,decorate, decoration={zigzag,amplitude=1.5pt,segment length = 3pt,pre length=2pt,post length=2pt}},
	rho/.style={dotted,semithick,shorten >=1pt,shorten <=1pt},
	renorm/.style={shape=circle,fill=white,inner sep=1pt},
	labl/.style={shape=rectangle,fill=white,inner sep=1pt},
	xi/.style={circle,fill=symbols!10,draw=symbols,inner sep=0pt,minimum size=1.2mm},
	xix/.style={crosscircle,fill=symbols!10,draw=symbols,inner sep=0pt,minimum size=1.2mm},
	xib/.style={circle,fill=symbols!10,draw=symbols,inner sep=0pt,minimum size=1.6mm},
	xibx/.style={crosscircle,fill=symbols!10,draw=symbols,inner sep=0pt,minimum size=1.6mm},
	not/.style={circle,fill=symbols,draw=symbols,inner sep=0pt,minimum size=0.5mm},
	>=stealth,
	}

\colorlet{symbols}{blue!90!black}

\makeatletter
\def\DeclareSymbol#1#2#3{\expandafter\gdef\csname MH@symb@#1\endcsname{\tikz[baseline=#2,scale=0.15]{#3}}%
\expandafter\gdef\csname MH@symb@#1s\endcsname{\scalebox{0.6}{\tikz[baseline=#2,scale=0.15]{#3}}}}
\def\<#1>{\csname MH@symb@#1\endcsname}
\makeatother

\DeclareSymbol{X}{-2.4}{\node[dot] {};}
\DeclareSymbol{1}{0}{\draw[white] (-.4,0) -- (.4,0); \draw (0,0)  -- (0,1.2) node[dot] {};}
\DeclareSymbol{2}{0}{\draw (-0.5,1.2) node[dot] {} -- (0,0.2) -- (0.5,1.2) node[dot] {};}
\DeclareSymbol{3}{0}{\draw (0,0) -- (0,1.2) node[dot] {}; \draw (-.7,1) node[dot] {} -- (0,0) -- (.7,1) node[dot] {};}
\DeclareSymbol{31}{-3}{\draw (0,0) -- (0,-0.7) node[root] {} -- (0.9,0.1) node[dot] {}; \draw (0,0) -- (0,1.2) node[dot] {}; \draw (-.7,1) node[dot] {} -- (0,0) -- (.7,1) node[dot] {};}
\DeclareSymbol{30}{-3}{\draw (0,0) -- (0,-0.7); \draw (0,0) -- (0,1.2) node[dot] {}; \draw (-.7,1) node[dot] {} -- (0,0) -- (.7,1) node[dot] {};}
\DeclareSymbol{10}{-3}{\draw (0,0) -- (0,-1); \draw (-.8,1) node[dot] {} -- (0,0);}
\DeclareSymbol{32}{-3}{\draw (0,0) -- (0,-0.7) -- (0.7,0) node[dot] {}; \draw  (0,-0.7)  -- (-0.7,0) node[dot] {}; \draw (0,0) -- (0,1) node[dot] {}; \draw (-.7,0.8) node[dot] {} -- (0,0) -- (.7,0.8) node[dot] {}; \node [root] at (0,-0.7) {};}
\DeclareSymbol{12}{-3}{\draw (0,0) -- (0,-1) -- (1,0) node[dot] {}; \draw (0,0) -- (0,-1) -- (-1,0) node[dot] {}; \draw (-.8,1) node[dot] {} -- (0,0);}
\DeclareSymbol{22}{-3}{\draw (0,0) -- (0,-0.7) -- (0.7,0) node[dot] {}; \draw (0,0) -- (0,-0.7) -- (-0.7,0) node[dot] {};\draw (-.5,0.8) node[dot] {} -- (0,0) -- (.5,0.8) node[dot] {};  \node [root] at (0,-0.7) {};}
\DeclareSymbol{20}{-3}{\draw (0,0) -- (0,-0.8);\draw (-.6,0.8) node[dot] {} -- (0,0) -- (.6,0.8) node[dot] {};}
\DeclareSymbol{32W}{-3}{\draw  [densely dotted, semithick] (-1.5,0)node[dot] {}  -- (0,-1.5) node[yy]{} -- (1.5,0) node[dot] {}; \draw (.8,-2.3)node{\tiny $y$} ;\draw[semithick] (0,0) -- (0,-1.5) node[yy] {}; \draw [densely dotted, semithick] (0,0) -- (0,1.8) node[dot] {}; \draw[densely dotted, semithick] (-1,1.5) node[dot] {} -- (0,0) -- (1,1.5) node[dot] {};}
\DeclareSymbol{32WW}{-3}{\draw  [densely dotted, semithick] (-1.5,0)node[dot] {}  -- (0,-1.5) node[yy]{} -- (1.5,0) node[dot] {}; \draw (.8,-2.3)node{\tiny $\bar{y}$} ;\draw[semithick] (0,0) -- (0,-1.5) node[yy] {}; \draw [densely dotted, semithick] (0,0) -- (0,1.8) node[dot] {}; \draw[densely dotted, semithick] (-1,1.5) node[dot] {} -- (0,0) -- (1,1.5) node[dot] {};}

\DeclareSymbol{s1}{0}{\draw[white] (-.4,0) -- (.4,0); \draw[symbols] (0,0)  -- (0,1.2) node[symb] {};}
\DeclareSymbol{s2}{0}{\draw[symbols] (-0.5,1.2) node[symb] {} -- (0,0) -- (0.5,1.2) node[symb] {};}
\DeclareSymbol{s3}{0}{\draw[symbols] (0,0) -- (0,1.2) node[symb] {}; \draw[symbols] (-.7,1) node[symb] {} -- (0,0) -- (.7,1) node[symb] {};}
\DeclareSymbol{s31}{-3}{\draw[symbols] (0,0) -- (0,-1) -- (1,0) node[symb] {}; \draw[symbols] (0,0) -- (0,1.2) node[symb] {}; \draw[symbols] (-.7,1) node[symb] {} -- (0,0) -- (.7,1) node[symb] {};}
\DeclareSymbol{s30}{-3}{\draw[symbols] (0,0) -- (0,-1); \draw[symbols] (0,0) -- (0,1.2) node[symb] {}; \draw[symbols] (-.7,1) node[symb] {} -- (0,0) -- (.7,1) node[symb] {};}
\DeclareSymbol{s32}{-3}{\draw[symbols] (0,0) -- (0,-1) -- (1,0) node[symb] {}; \draw[symbols] (0,0) -- (0,-1) -- (-1,0) node[symb] {}; \draw[symbols] (0,0) -- (0,1.2) node[symb] {}; \draw[symbols] (-.7,1) node[symb] {} -- (0,0) -- (.7,1) node[symb] {};}
\DeclareSymbol{s22}{-3}{\draw[symbols] (0,0.3) -- (0,-1) -- (1,0) node[symb] {}; \draw[symbols] (0,0.3) -- (0,-1) -- (-1,0) node[symb] {};\draw[symbols] (-.7,1) node[symb] {} -- (0,0.3) -- (.7,1) node[symb] {};}
\DeclareSymbol{s20}{-3}{\draw[symbols] (0,0) -- (0,-1);\draw[symbols] (-.7,1) node[symb] {} -- (0,0) -- (.7,1) node[symb] {};}
\DeclareSymbol{s21}{-3}{\draw[symbols] (0,0) -- (0,-1);\draw[symbols] (-.7,1) node[symb] {} -- (0,0) -- (.7,1) node[symb] {}; \draw[symbols] (0,0) -- (0,-1) -- (1,0) node[symb] {};}
\DeclareSymbol{s12}{-3}{\draw[symbols] (0,0) -- (0,-1) -- (1,0) node[symb] {}; \draw[symbols] (0,0) -- (0,-1) -- (-1,0) node[symb] {}; \draw (-.8,1) node[symb] {} -- (0,0);}
\DeclareSymbol{s11}{-3}{\draw[symbols] (0,0) -- (0,-1) -- (-1,0) node[symb] {}; \draw (-.8,1) node[symb] {} -- (0,0);}
\DeclareSymbol{s10}{-3}{\draw[symbols] (0,0) -- (0,-1); \draw[symbols] (-.8,1) node[symb] {} -- (0,0);}

\definecolor{darkergreen}{rgb}{0.0, 0.5, 0.0}

\numberwithin{equation}{section}
\def\theequation{\arabic{section}.\arabic{equation}}
\newcommand{\be}{\begin{eqnarray}}
\newcommand{\ee}{\end{eqnarray}}
\newcommand{\ce}{\begin{eqnarray*}}
\newcommand{\de}{\end{eqnarray*}}
\newtheorem{theorem}{Theorem}[section]
\newtheorem{lemma}[theorem]{Lemma}
\newtheorem{remark}[theorem]{Remark}
\newtheorem{definition}[theorem]{Definition}
\newtheorem{proposition}[theorem]{Proposition}
\newtheorem{Examples}[theorem]{Example}
\newtheorem{corollary}[theorem]{Corollary}

\newcommand{\LL}{\mathscr{L}}
\newcommand{\UU}{\mathscr{U}}

\def\Wick#1{\,\colon\!\! #1 \!\colon}

\def\PPhi{\mathbf{\Phi}}

\def\eps{\varepsilon}

\def\p{\partial}

\def\l{\lambda}
\def\la{\langle}
\def\ra{\rangle}
\def\[{{\Big[}}
\def\]{{\Big]}}
\def\({{\Big(}}
\def\){{\Big)}}

\def\bx{{\mathbf{x}}}
\def\tr{\mathrm {tr}}

\def\dif{{\mathord{{\rm d}}}}

\def\min{{\mathord{{\rm min}}}}

\def\no{\nonumber}
\def\={&\!\!=\!\!&}

\def\bC{{\mathbf C}}

\def\cI{{\mathcal I}}

\def\cZ{{\mathcal Z}}

\def\1{{\mathbf{1}}}

\def\sD{{\mathscr D}}

\def\E{\mathbf E}

\def\geq{\geqslant}
\def\leq{\leqslant}
\def\ge{\geqslant}
\def\le{\leqslant}

\def\eps{\varepsilon}

\def\p{\partial}

\def\l{\lambda}
\def\la{\langle}
\def\ra{\rangle}
\def\[{{\Big[}}
\def\]{{\Big]}}
\def\({{\Big(}}
\def\){{\Big)}}

\def\bx{{\mathbf{x}}}
\def\tr{\mathrm {tr}}

\def\dif{{\mathord{{\rm d}}}}

\def\min{{\mathord{{\rm min}}}}

\def\no{\nonumber}
\def\={&\!\!=\!\!&}
\def\bt{\begin{theorem}}
\def\et{\end{theorem}}
\def\bl{\begin{lemma}}
\def\el{\end{lemma}}
\def\br{\begin{remark}}
\def\er{\end{remark}}
\def\bx{\begin{Examples}}
\def\ex{\end{Examples}}
\def\bd{\begin{definition}}
\def\ed{\end{definition}}
\def\bp{\begin{proposition}}
\def\ep{\end{proposition}}
\def\bc{\begin{corollary}}
\def\ec{\end{corollary}}

\def\geq{\geqslant}
\def\leq{\leqslant}
\def\ge{\geqslant}
\def\le{\leqslant}

 \def\R{\mathbb R}
 \def\R{\mathbb R}    
\def\N{\mathbb N}

\allowdisplaybreaks

\begin{document}

\title[An SPDE approach to large N problems  in QFT: a survey]{A stochastic PDE approach to large N problems  in quantum field theory: a survey}

\author[Hao Shen]{Hao Shen}
\address[H. Shen]{Department of Mathematics, University of Wisconsin - Madison, USA}
\email{pkushenhao@gmail.com}
%

\begin{abstract}
In this survey we review some recent rigorous results
on large N problems in quantum field theory, 
stochastic quantization and singular stochastic PDEs, and
their mean field limit problems.
In particular we discuss the $O(N)$ linear sigma model on two and three dimensional torus. The stochastic quantization procedure leads to a coupled system of $N$ interacting $\Phi^4$ equations. 
In $d=2$, we show uniform in $N$ bounds for the dynamics and convergence to a mean-field singular SPDE.  
For large enough mass or small enough coupling, 
the invariant measures (i.e. the $O(N)$ linear sigma model)
 converge to the massive Gaussian free field, the unique invariant measure of the mean-field dynamics, 
 in a Wasserstein distance.  
 We also obtain tightness for certain $O(N)$ invariant observables as random fields in suitable Besov spaces as $N\to \infty$, along with exact descriptions of the limiting correlations.
In $d=3$, the estimates become more involved since the equation is more singular. 
We discuss in this case how to prove convergence to the massive Gaussian free field. 
%
The proofs of these results build on the recent progress of singular SPDE theory
and combine many new techniques such as uniform in $N$ estimates and dynamical mean field theory.
These are based on joint papers with Scott Smith, Rongchan Zhu and Xiangchan Zhu.
\end{abstract}


\date{\today}

\maketitle


\section{Introduction and physics background}

In quantum field theory (QFT), there are important models where the ``field'' takes values in an $N$ dimensional space.
A natural question is that what is the asymptotic behavior of such models as dimensionality of the target space $N\to \infty$.
The perspective of investigating these questions in this article is called stochastic quantization. 
Stochastic quantization  is a procedure which formulates a quantum field theory in terms of a stochastic PDE (SPDE). 
More precisely, for a quantum field theory where the target space has dimension $N$, its stochastic quantization yields
a system of $N$ coupled SPDEs.
This procedure then brings in many powerful analytical tools, including the recent developments of SPDE solution theories and their a priori estimates, and mean field theory techniques. 
The purpose of this article is to survey the recent results along this direction.

For most of the article, we will focus on the linear  $\sigma$-model in $2$ and $3$ dimensions. But we will first briefly review the origin of the large N problems in physics.

\subsection{Large N problem in QFT}
\label{sec:LargeN}

Large $N$ methods (or ``$1/N$ expansions'') in theoretical physics are ubiquitous and
are generally applied to models where dimensionality of the target space is large.
It was first noted by Stanley \cite{stanley1968spherical} for lattice spin models (rather than QFT),
that they exhibit considerable simplification as $N$ (the number of spin components) becomes large: 
at least at the level of free energy, as $N\to \infty$, these models become the ``spherical model'', which is a solvable model
 introduced by Berlin and Kac in 1952 \cite{berlin1952spherical}. See Baxter's book 
 \cite{baxter2016exactly}, Section 5, for an expository discussion on spherical model.
 Later, the large $N$ convergence to the spherical model
 was established at the level of free energies \cite{kac1971spherical}
 and then of correlations  \cite{MR471851,MR982418}.
Again in the context of spin systems,
 \cite{brezin1973critical} found a systematic way to expand in powers of $1/N$, 
and obtained $1/N$ order
corrections 
 to the spherical model.

In 1973 Wilson exploited 
this idea  in quantum field theory \cite{wilson1973quantum}.
In particular Wilson  
studied both bosonic (linear $\sigma$-model) and fermionic models.
The linear $\sigma$-model
is the
 $N$-component generalization of the $\Phi^4_d$ model, given by the  (formal) measure 
\begin{equation}\label{e:Phi_i-measure}
\dif\nu^N(\Phi)\eqdef \frac{1}{Z_N}\exp\bigg(-\int_{\mathbb T^d} \sum_{j=1}^N|\nabla \Phi_j|^2+m \sum_{j=1}^N\Phi_j^2
+\frac{\lambda}{2N} \Big(\sum_{j=1}^N\Phi_j^2\Big)^2 \dif x\bigg)\mathcal D \Phi
\end{equation}
over $\R^N$ valued fields $\Phi=(\Phi_1,\Phi_2,...,\Phi_N)$,
where $Z_N$ is a normalization constant and $\mathcal D \Phi$ is the formal ``Lebesgue measure''.
We will discuss below the physics predictions  and how they are related with our rigorous results.  Soon later, Coleman--Jackiw--Politzer \cite{coleman1974}
calculated the effective potential and studied spontaneous symmetry breaking of this model for $d\le 4$, see also \cite{abbott1976bound} for more extensive analysis in $d=4$. 
\footnote{Remark that in this survey we will not discuss phase transition or 
spontaneous symmetry breaking, although it is an interesting topic.}

Around the same time, Gross and Neveu \cite{gross1974dynamical} studied the $1/N$ expansion for a two-dimensional
fermionic model, where $N$ Dirac fermions interact 
via a quartic term 
(``four-fermion interactions''): namely the Lagrangian is given by 
$\sum_{i=1}^N\bar\psi_i (i\slashed{\partial})\psi_i + \frac{\lambda}{N} (\sum_{i=1}^N\bar\psi_i \psi_i)^2 $.
This model had also been touched in Wilson's paper \cite{wilson1973quantum}. It was observed that the model has asymptotic freedom,
symmetry breaking (for a discrete chiral symmetry) and produces a fermion mass.

Another class of models considered from the 70s
is the nonlinear $\sigma$-models in two dimensions.
For instance, 
 see \cite{bardeen1976phase} which considered
models where $\Phi$ takes values in a sphere in $\R^N$.
It can be viewed as 
the limit of  $\R^N$ valued $\sigma$-model  
with a potential $K (\Phi\cdot\Phi - a^2)^2$
in a formal limit $K\to\infty$.
The idea was also extended to $CP^{N-1}$ type nonlinear $\sigma$-models \cite{d1979confinement,d19781n} (where more general homogenous space $G/H$ valued models were also formulated).

Matrix-valued field theories are considered to be more difficult. In 1974 t'Hooft
\cite{tHooft1974planar} gave the solution of the two dimensional Yang-Mills model with $SU(N)$ gauge group,
as $N\to \infty$.
After this seminal work by t'Hooft, the large $N$ behavior of gauge theories was extensively studied by physicists,
see e.g. \cite{makeenko1979exact,polyakov1980gauge,kazakov1980non,kazakov1981wilson}.
See also \cite{gopakumar1995mastering} (and the references therein) for a broad class of matrix models.

It is certainly impossible to list all the vast literature,
instead we refer to  \cite{LargeN1993} for an edited comprehensive collection of articles on large
$N$ as applied to a wide spectrum of problems in quantum field theory and statistical
mechanics; see also the review articles 
\cite{witten1980} and \cite{moshe2003} 
or Coleman's excellent book  \cite[Chapter~8]{coleman1988} for summaries of the progress.

\subsection{Perturbation theory heuristics}
\label{sec:Perturbation}

A standard approach in quantum field theory is perturbation theory, namely, calculating physical quantities of a given model
by perturbative expansion in a coupling constant of the model.
In large $N$ problems, terms  in perturbation theory can have different orders in $1/N$,
and the ``$1/N$ expansion'' is a re-organization of the series in the parameter $1/N$,
 with each term typically being a (formal) sum of infinitely many orders of the ordinary perturbation theory. 
 
To explain this formalism a bit more, we first recall the usual perturbation
theory for $\Phi^4$ model, with $N=1$ and $m=0$ (for simplicity) in  \eqref{e:Phi_i-measure}. This could be found in a standard QFT textbook, e.g.  Peskin and  Schroeder \cite{MR1402248}. Loosely speaking, the perturbative calculation of for instance a two-point correlation of $\Phi$ is given by a sum of Feynman graphs with two external legs and degree-$4$ internal vertices that represents the quartic interaction $\Phi^4$, such as
 
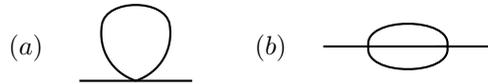
\begin{center}
$(a)\quad$
\begin{tikzpicture}[scale=1.5,baseline=10]
\draw[thick] (0,0) to  (1,0);
\draw[thick] (0.5,0)  [bend left=40] to (0.2,0.5);
\draw[thick]  (0.2,0.5)  [bend left=80] to (0.8,0.5);
\draw[thick]  (0.8,0.5) [bend left=40] to (0.5,0);
\end{tikzpicture}
$\qquad (b)\quad$
\begin{tikzpicture}[scale=1.5,baseline=-3]
\draw[thick] (0,0) to  (0.4,0);
\draw[thick] (0.4,0)  [bend left=100] to (1.1,0);
\draw[thick] (0.4,0) to (1.1,0);
\draw[thick] (0.4,0)  [bend right=100] to (1.1,0);
\draw[thick] (1,0) to  (1.5,0);
\end{tikzpicture}
\end{center}
For the vector model \eqref{e:Phi_i-measure} with $N>1$ which was studied
by Wilson \cite{wilson1973quantum}, the  perturbative calculation of the two-point correlation of $\Phi_i$ (with $i$ fixed) is also given by the sum of Feynman graphs with two external legs and degree-$4$ internal vertices, but now each vertex carrying two distinct summation variables
and a factor $1/N$ that represents the interaction $\frac1N \sum_{j,k}\Phi_j^2 \Phi_k^2$ term in  \eqref{e:Phi_i-measure}, such as
(a) (b) below

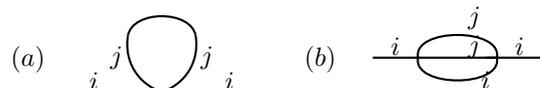
\begin{center}
$(a)\quad$
\begin{tikzpicture}[scale=1.5,baseline=10]
\draw[thick] (0,0) to  (1,0);
\draw[thick] (0.5,0)  [bend left=40] to (0.2,0.5);
\draw[thick]  (0.2,0.5)  [bend left=80] to (0.8,0.5);
\draw[thick]  (0.8,0.5) [bend left=40] to (0.5,0);
\node at (-0.1,0.1) {$i$}; \node at (1.1,0.1) {$i$};
\node at (0.1,0.3) {$j$}; \node at (0.9,0.3) {$j$};
\end{tikzpicture}
$\qquad (b)\quad$
\begin{tikzpicture}[scale=1.5,baseline=-3]
\draw[thick] (0,0) to  (0.4,0);
\draw[thick] (0.4,0)  [bend left=100] to (1.1,0);
\draw[thick] (0.4,0) to (1.1,0);
\draw[thick] (0.4,0)  [bend right=100] to (1.1,0);
\draw[thick] (1,0) to  (1.5,0);
\node at (0.2,0.1) {$i$}; \node at (1.3,0.1) {$i$};
\node at (0.9,0.1) {$j$}; \node at (0.9,0.35) {$j$}; \node at (1,-0.2) {$i$};
\end{tikzpicture}
\end{center}
Heuristically,
 graph (a) is of order $\frac1N \sum_j \approx O(1)$ and
  graph (b) is of order $\frac{1}{N^2} \sum_j \approx O(\frac1N)$.
 The philosophy of  \cite{wilson1973quantum} is that graphs with ``self-loops'' such as (a) get {\it cancelled} by Wick renormalization,
 \footnote{We will have more precise discussion on Wick renormalization in Section~\ref{sec:Phi42}.}
 and all other graphs with  internal vertices  such as (b) are at least of order $O(1/N)$ and thus vanish,
so the theory would be asymptotically Gaussian free field.
(We rigorously prove this in Section~\ref{sec:GFF2d} and Section~\ref{sec:GFF3d} in 2D and 3D.)

Since the model \eqref{e:Phi_i-measure}
is invariant under any $O(N)$ rotation of $\Phi$, 
it is more natural to study the $O(N)$-invariant observables
than the field itself, such as
$\sum_{j=1}^N\Phi_j^2$.
One needs to find a suitable scaling for such observables
to see interesting fluctuations.
In the ``law of large numbers behavior'' scaling $\frac1N \sum_{j=1}^N\Phi_j^2$
one would expect a deterministic limit.
Interesting fluctuation then arises in the next scale.

Following the above heuristic, we assume that the model is Wick renormalized
and so the perturbation series does not contain self-loops, and 
we consider
 $\frac{1}{\sqrt N} \sum_{j=1}^N \Wick{\Phi_j^2}$.
Note that here  $\Wick{\Phi_j^2}$ is the Wick square.
The bubbled diagrams of the following form
are all the $O(1)$ contributions to its two-point correlation:
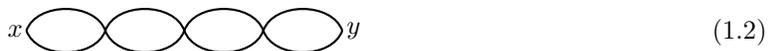
\begin{equ}[e:bubbles]
\begin{tikzpicture}[scale=1.5,baseline=-3]
\draw[thick] (0,0)  [bend left=70] to (0.7,0);
\draw[thick] (0,0)  [bend right=70] to (0.7,0);
\draw[thick] (0.7,0)  [bend left=70] to (1.4,0);
\draw[thick] (0.7,0)  [bend right=70] to (1.4,0);
\draw[thick] (1.4,0)  [bend left=70] to (2.1,0);
\draw[thick] (1.4,0)  [bend right=70] to (2.1,0);
\draw[thick] (2.1,0)  [bend left=70] to (2.8,0);
\draw[thick] (2.1,0)  [bend right=70] to (2.8,0);
\node at (-0.1,0) {$x$};\node at (2.9,0) {$y$};
\end{tikzpicture}
\end{equ}
In fact, it is known in physics literature that these ``bubble chains''
are dominant (after Wick renormalization) contributions for this model:
for instance if we consider another observable
$ \frac{1}{N} \Wick{\big(\sum_{i=1}^N\Phi_i^2\big)^2}$,
the leading contribution to its one-point correlation
would then be the ``bubble loops'' i.e. the above type 
of diagrams with $x,y$ identified.

\subsection{Dyson-Schwinger and dual field method}
\label{sec:DS}

Besides directly examining  the perturbation theory, alternative (and more systematic) methodologies  of analyzing such expansion were discovered in physics. We briefly recall some of them.

In \cite{symanzik1977} Schwinger-Dyson equations were used to study
models of the type \eqref{e:Phi_i-measure}. These are essentially a type of integration by parts formulae. 
They yield a hierarchy of relations for correlation functions of different orders.
This hierarchy does not close, so it is somewhat similar with perturbation theory;
but one advantage of using Schwinger-Dyson equations 
is that one could write a truncation remainder in terms of
again correlations of the field. This will be useful for rigorous studies (see \cite{MR578040}), and we will come back to this point in Section \ref{sec:observables}.

Schwinger-Dyson equations are particularly powerful when we
move from vector  valued models of the type  \eqref{e:Phi_i-measure}
to matrix valued models, such as Yang--Mills.
In the context of Yang-Mills,
\cite{makeenko1979exact} discovered 
relations between Wilson loop observables among different loops,
now called Makeenko--Migdal  equations or master loop equations.
These are also broadly used in the study of random matrices,
for instance see the book \cite{guionnet2019asymptotics} by Guionnet.

The problem of the  perturbative argument discussed in Sec.~\ref{sec:Perturbation}
(besides being non-rigorous)
is that keeping track of all the diagrams
and proving that there is no other diagrams appearing
at the leading order is a combinatoric challenge;
this challenge also exists when iteratively applying Schwinger-Dyson equations.
A better way is to introduce a ``dual'' field or ``auxiliary'' field,
which we briefly review now.
A very nice exposition of this idea can be found in 
Coleman's book  \cite[Chapter~8]{coleman1988},
(see also the original paper \cite{coleman1974}, or the review \cite[Section~2]{moshe2003}).
The exact ways of implementing this idea vary among the literature but roughly speaking 
one will shift the action in \eqref{e:Phi_i-measure}
 by a quadratic term involving a new field $\sigma$ when calculating
 an expectation w.r.t. the measure $\dif\nu^N$:
 \begin{equ}
\int (\cdots ) \dif\nu^N(\Phi)
\to
\int \exp\bigg( -\int  \frac{N}{2\lambda} (\sigma - \frac{\lambda}{N}\sum_{j=1}^N\Phi_j^2)^2\bigg)
\dif\sigma \dif\nu^N(\Phi)\;.
\end{equ}
This is because the integral over $\sigma$
is simply a Gaussian integral which yields 
a constant that can be absorbed into the partition function.
We then obtain $\exp(-\mathcal S)$ where
$$
\mathcal S (\sigma,\Phi)=
\int \sum_{j=1}^N|\nabla \Phi_j|^2+m \sum_{j=1}^N\Phi_j^2
+ \frac{N}{2\lambda} \sigma^2 - \frac12  \sum_{j=1}^N \sigma \Phi_j^2 \; \dif x
$$
in particular the quartic term is cancelled. 
Still considering two-point correlation of $\Phi_i$ for instance,
one has the following type of diagrams
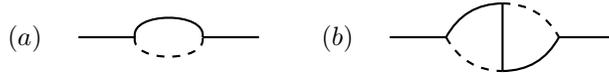
\begin{center}
$(a)\quad$
\begin{tikzpicture}[scale=1.5,baseline=-3]
\draw[thick] (-0.1,0) to  (0.4,0);
\draw[thick] (0.4,0)  [bend left=100] to (1,0);
\draw[thick,dashed] (0.4,0)  [bend right=100] to (1,0);
\draw[thick] (1,0) to  (1.5,0);
\end{tikzpicture}
$\qquad(b)\quad$
\begin{tikzpicture}[scale=1.5,baseline=-3]
\draw[thick] (0,0) to  (0.5,0);
\draw[thick] (0.5,0)  [bend left=30] to (1,0.3);
\draw[thick,dashed] (1,0.3)  [bend left=30] to (1.5,0);
\draw[thick,dashed] (0.5,0)  [bend right=30] to (1, -0.3);
\draw[thick] (1, -0.3)  [bend right=30] to (1.5,0);
\draw[thick] (1,0.3) to  (1,-.3);
\draw[thick] (1.5,0) to  (2,0);
\end{tikzpicture}
\end{center}
Here the dashed line 
represent the covariance of the ``Gaussian measure''
$\exp(-\int \frac{N}{2\lambda} \sigma(x)^2 \,\dif x)$,
which is simply $\frac{\lambda}{N} \delta$.
Each internal vertex is degree 3,
because of the $\sigma \Phi_j^2$ interaction term.

One could study this theory of $(\sigma,\Phi)$ but
as explained in
  \cite{coleman1988}, 
 since the above action is quadratic in $\Phi$, one could ``integrate out''
  $(\Phi_1,\cdots,\Phi_N)$
 and obtain an ``effective action'' for $\sigma$. 
 This Gaussian integral would however yields
 a non-local interaction because of the term $|\nabla\Phi_j|^2$,
 namely a determinant $\det(-\Delta+m -\frac{\sigma}2)^{N/2}$,
 or $\det(1 -\frac12 (-\Delta+m)^{-1} \sigma)^{N/2}$ if we factor out the determinant of the massive Laplacian.
Many physical properties can be deduced from 
this dual field theory for $\sigma$,
for instance one could exponentiate this determinant 
to obtain a measure of the form $e^{N F(\sigma)}$ where $F$ has the form 
$\tr \log (...)$ 
so that large $N$ limit becomes a question of critical point analysis for $F$.
The $1/N$ expansion for correlations will also turn out to be a simple perturbative expansion in the dual picture.

\subsection{Some previous rigorous results on large $N$}

Rigorous  study of large $N$ in mathematical physics was initiated  by Kupiainen \cite{MR574175,MR578040} around 1980. See also his review \cite{MR582622}.
The literature most related to the present article is \cite{MR578040}, which studied the linear sigma model in continuum in $d=2$
(with Wick renormalization)
and proved that the $1/N $ expansion of the  pressure (i.e. vacuum energy, or log of partition per area)  is asymptotic,
 and each order in this expansion
  can be described by sums of infinitely many Feynman diagrams of certain types.
 This paper uses both Dyson-Schwinger equations and dual field method,
together with constructive field theory methods such as chessboard estimates 
\cite{MR436829} based on reflection positivity.
 Borel summability of $1/N $ expansion of Schwinger functions for this model was discussed in \cite{MR661137}.
 
Most of the other rigorous results are proved on lattices and thus ignore 
the small scale singularity (which is one of the key difficulties in our work).
We also review some of these results here even if they are not necessarily directly related with our results.
In \cite{MR574175}, Kupiainen  
considered the $N$-component
nonlinear sigma model (fields take values in a sphere in $\R^N$)
on the $d$ dimensional lattice with fixed lattice spacing.
As mentioned above this simplifies to ``spherical model'' as $N\to \infty$,
and for the spherical model it is known that 
the critical temperature is $0$ when $d\le 2$
and is strictly positive when $d>2$.
\cite{MR574175}
proved that  the large $N$ expansion of correlation functions and free energy 
 is asymptotic above the spherical model critical temperature (meaning all temperature if $d\le 2$),
  and mass gap was established for these temperatures when $N$ is sufficiently large. The main idea therein was the dual field representation
(for which the heuristic was discussed in Sec.~\ref{sec:DS})
which yields nonlocal expectations, and this technical difficulty was solved using random walk expansion ideas from Brydges--Federbush \cite{MR496278}.
Asymptoticity was later extended to Borel summability by \cite{MR678004}.

In the 90s, Kopper--Magnen--Rivasseau \cite{MR1328264}
proved the mass gap for fermions in the Gross--Neveu model when $N$ is sufficiently large, and the model has at least two pure phases. They again work with the dual field as illustrated above
(which ends up with a determinant of Dirac operator) with suitable ultraviolet regularization and apply cluster and Mayer expansion techniques. Note that
the Gross--Neveu model action has a (discrete) chiral symmetry 
which prohibits a mass but the mass is acquired in a symmetry breaking.
\cite{MR1686539} (on lattice $\Z^2$)
\cite{MR1686543} (on $\R^2$ but with certain ultraviolet cutoff and a ``soft constraint'' $K (\Phi\cdot\Phi - 1)^2$ in the potential) and  then revisited 
the nonlinear $\sigma$ model and showed mass gap for $N$ sufficiently large,
and they both work with dual fields and determinants
using cluster expansions.

Large $N$ limit and expansion for Yang-Mills model
has drawn much attention, following t'Hooft's 1974 work
\cite{tHooft1974planar}.
In two space dimensions (continuum plane), L{\'e}vy
\cite{Levy11}  
proved 
convergence of Wilson loop observables to master field in the large $N$ limit, with orthogonal, unitary and symplectic structure groups.
The master field is a deterministic object,
which was described conjecturally by I.M. Singer in 1995
\cite{MR1373007} who pointed out its connections with free probability.
L{\'e}vy's work
\cite{Levy11} also rigorously establishes integration by parts i.e. Schwinger-Dyson equations
 both for finite $N$ and in the large $N$ limit  for the expectations of Wilson loops, which in this context are called the Makeenko-Migdal equations.
 Note that two dimensional Yang-Mills model has special integrability structures,
which allows one to make sense of the random Wilson loop observables in continuum  - this was established by L{\'e}vy in \cite{MR2006374,MR2667871} which \cite{Levy11} relies on.
Around the same time Anshelevich-Sengupta \cite{MR2864481} also gave a construction of the master field  in the large $N$ limit,
via a different
approach which is based on the use of free white noise and of free stochastic calculus; these methods were developed earlier by Sengupta in \cite{MR1346931,MR2757706}.
On the lattice of any dimension,
we refer to 
 \cite{MR3919447} (resp. \cite{chatterjee2016}) for
computation of correlations of the Wilson loops 
 in the large $N$ limit (resp. $1/N$ expansion) which relates to  string theory.

%
%
%

We remark that while the  
Schwinger-Dyson equations for vector models (as we will discuss more 
in Section \ref{sec:observables}) are easy to derive,
Makeenko--Migdal equations for Yang-Mills models are more difficult to derive (rigorously).
The first rigorous version was established for two-dimensional Yang-Mills model in continuum in \cite{Levy11}, 
and alternative proofs and extensions were given in \cite{MR3554890,MR3982691,MR3631396,MR3613519} on plane or surfaces. 
On lattice, for finite $N$ and in any dimension, these were  derived by 
Chatterjee \cite{MR3919447} and  \cite{Jafar},  and then in
\cite{SSZloop} using stochastic ODEs on Lie groups.

%

\section{Stochastic quantization and SPDE}\label{sec:SQ}

Stochastic quantization
 refers a procedure to turn Euclidean quantum field theories to 
 singular SPDEs, and thus bring stochastic analysis and PDE methods into the study. 
 They were introduced by Parisi and Wu in \cite{ParisiWu}. Given an action $\mathcal S(\Phi)$ which is a functional of $\Phi$,
one considers a gradient flow of $\mathcal S(\Phi)$  perturbed by space-time white noise $\xi$:
\begin{equ}[e:SQE]
\partial_t \Phi = - \frac{\delta \mathcal S(\Phi)}{\delta \Phi} + \xi \;.
\end{equ}
Here $\frac{\delta \mathcal S(\Phi)}{\delta \Phi}$ is  the variational derivative of the functional $\mathcal S(\Phi)$;
for instance, when $\mathcal S(\Phi)=\frac12\int (\nabla \Phi)^2 dx$
is the Dirichlet form, $\frac{\delta \mathcal S(\Phi)}{\delta \Phi} =-\Delta \Phi$
and \eqref{e:SQE} boils down to the stochastic heat equation
$$
  \partial_t \Phi = \Delta \Phi   + \xi\;.
$$
Note that 
 $\Phi$ can be also multi-component fields, with $\xi$ being likewise multi-component. 
Another  example is  the ``dynamical $\Phi^4$ equation''
\begin{equ}[e:Phi4-formal]
  \partial_t \Phi = \Delta \Phi - \lambda\Phi^3  + \xi
\end{equ}
and it arises from the gradient of $\int \frac12 |\nabla\Phi|^2+ \frac{\lambda}4 \Phi^4\,dx$, that is, the model \eqref{e:Phi_i-measure} with $m=0,N=1$.

The significance of these ``stochastic quantization equations'' \eqref{e:SQE} is that
given an action $\mathcal S(\Phi)$,
the formal measure which defines an Euclidean quantum field theory
\begin{equ}[e:QFT]
\tfrac{1}{Z}e^{-\mathcal S(\Phi)} D\Phi
\end{equ}
 is formally an invariant measure\footnote{One even simpler example in stochastic ordinary differential equations is given by the Ornstein-Uhlenbeck process $dX_t=-\frac12 X_tdt+dB_t$ where $B_t$ is the Brownian motion, and its invariant measure is the (one-dimensional) Gaussian measure $\frac{1}{\sqrt{2\pi}}e^{-\frac{X^2}{2}}dX$.}
 for Eq. \eqref{e:SQE}.
 Here  $ D\Phi$ is the formal Lebesgue measure and
 $Z$ is a ``normalization constant''.
We emphasize that \eqref{e:QFT} are only formal measures because, among several other reasons, there is no ``Lebesgue measure'' $D\Phi$ on an  infinite dimensional space and it is a priori not clear at all if the measure can be normalized.
The task of {\it constructive quantum field theory} is to give precise meaning or constructions to these formal measures. This area of study has yielded numerous deep results which we do not mean to give a full list of literature; we only refer to 
the book \cite{MR887102,MR1773042}.

Regarding  the well-posedness of \eqref{e:Phi4-formal}, in two dimensions, two classical works are \cite{MR1113223} where martingale solutions are constructed and \cite{DD03} where local strong solutions are obtained, as well as the more recent approach to global well-posedness in \cite{MW17}. 
In three dimension the local strong solution is obtained as an application of regularity structures \cite{Hairer14} and then \cite{MR3846835}
as an application of  paracontrolled analysis.

Given the very recent progress of SPDE,
a new  approach to construct 
the measure of the form \eqref{e:QFT} is to construct the {\it long-time} solution to the stochastic PDE \eqref{e:SQE} and average the distribution of the solution over time.
This approach has shown to be successful for the $\Phi^4$ model in $d\le 3$ 
in a series of very recent works, which starts with
 \cite{MW18} on the 3D torus  where a priori estimates were obtained to rule out the possibility of finite time blow-up.
Then in \cite{GHglobal,GH21} established a priori estimates for solutions on the full space $\R^3$ yielding the construction of $\Phi^4$ quantum field theory 
on the whole $\R^3$, as well as verification of some key properties that this invariant measure must satisfy as desired by physicists such as reflection positivity. See also \cite{AK17}. Similar uniform a priori estimates are obtained by \cite{moinat2020space} using  maximum principle.

  Large $N$ problems in the stochastic quantization formalism
 have also been discussed (on a heuristic level) in the physics literature, for instance  \cite{Alfaro1983,AlfaroSakita},
 \cite[Section~8]{damgaard1987}.
 The review on large $N$ by
Moshe and Zinn-Justin
 \cite[Section~5.1]{moshe2003}  is close to our setting; it makes an ``ansatz'' that $\frac{1}{N}\sum_{j=1}^N\Phi_j^2$ in the equation would self-average in the large $N$ limit to a constant; our results justified this ansatz and in the non-equilibrium setting generalizes it.

\section{Mean field limits for coupled stochastic dynamics}
\label{sec:SDE}

As we will see, 
 the methods in the proofs of our main theorems borrow some ingredients from
mean field limit theory (MFT).
To the best of our knowledge, the study of mean field problems originated from McKean \cite{MR0233437}.
Typically, a mean field problem is concerned with a system of $N$ particles interacting with each other,
which is often modeled by a system of stochastic {\it ordinary} differential equations, for instance, driven by independent Brownian motions. A prototype of such systems has the form $\dif X_i = \frac{1}{N} \sum_j f(X_i,X_j) \dif t + \dif B_i$, see for instance the classical reference by Sznitman \cite[Sec~I(1)]{MR1108185},
 and in the $N\to \infty$ limit one could obtain decoupled SDEs each interacting with the law of itself:  we will briefly review this calculation below.

In simple situations the interaction $f$ is assumed to be ``nice'', for instance globally Lipschitz (\cite{MR0233437}); much of the literature aims to prove such limits under more general assumptions on the interaction, see  \cite{MR1108185} for a survey.
\footnote{In the context of SDE systems, one also considers the empirical measures of the particle configurations, and aims to show their convergence as  $N\to \infty$ to the McKean--Vlasov PDEs, which are typically deterministic. 
Note that in this paper we do not consider the ``analogue'' of McKean-Vlasov PDE (which would be infinite dimensional) in the context of our model.}

We note that mean field limits are studied under much broader frameworks or scopes of applications,
such as mean field limit in the context of rough paths  (e.g. \cite{MR3299600,BCD18,CDFM18}),
mean field games (e.g. survey \cite{MR2295621}), quantum dynamics (e.g. \cite{MR2680421} and references therein).
We do not intend to have a comprehensive list, but rather refer to survey articles \cite{MR3468297,MR3317577} and the book \cite[Chapter~8]{spohn2012large}
besides \cite{MR1108185}.

The study of mean field limit for SPDE systems also has precursors, see for instance the book \cite[Chapter~9]{MR1465436} or \cite{MR3160067}. However these results make strong assumptions on the interactions of the SPDE systems such as linear growth and globally Lipschitz drift, and certainly do not cover the singular regime where renormalization is required as in our case.

Before we delve into the mean field limit for singular SPDE systems,
we briefly review a very simple example regarding SDEs, which is taken from 
Sznitman \cite[Sec~I(1)]{MR1108185}. We will see that some of the simple ideas will be reflected when we move on to the much more technical SPDE setting.
Consider
$$
dX_i = \frac{1}{N} \sum_{j=1}^N f(X_i,X_j) dt+ dB_i \;, \qquad (i=1,...,N)
$$
where $B_i$ ($i=1,...,N$) are independent Brownian motions. One can show that as $N\to \infty$, the solution converges to that of (``McKean--Vlasov limit'')
\begin{equ}[e:barXSDE]
d \bar X_i = \int_{\mathbf R} f(\bar X_i, y) u_t(y)dydt
+ dB_i \;, \qquad (i=1,...,N) 
\end{equ}
where $u_t = \mbox{Law} (\bar X_i(t))$. 
 Note that this is a decoupled system now.
To prove the above result, we subtract the two systems and get
\begin{equs}[e:sznit]
\partial_t (X_i & -\bar X_i) 
= \frac{1}{N} \sum_{j=1}^N f(X_i,X_j) 
- \int_{\mathbf R} f(\bar X_i, y) u_t(y) dy
\\
&-\frac{1}{N} \sum_{j=1}^N f(\bar X_i,X_j) 
+\frac{1}{N} \sum_{j=1}^N f(\bar X_i,X_j) 
-\frac{1}{N} \sum_{j=1}^N f(\bar X_i,\bar X_j) 
+\frac{1}{N} \sum_{j=1}^N f(\bar X_i,\bar X_j) \;.
\end{equs}
Here the terms in the second line 
are put ``by hand'' and they obviously sum to zero.
We then compare the 1st and 3rd terms on the RHS,
and the 4th  and  5th terms, and the 2nd and the last term;
we show that the differences all vanish as $N\to \infty$.
The most important case is the following:
\begin{equ}[e:1N]
\E \Big(
\frac1N \sum_{j=1}^N 
\Big( f(\bar X_i,\bar X_j) 
 - \int_{\mathbf R} f(\bar X_i, y) u_t(y) dy \Big)
 \Big)^2 \sim 1/N\;.
\end{equ}
Note that the nontriviality here is that 
the left hand side ``appears to be'' $1/N^2$ times two summations,
which would be $O(1)$ as $N\to \infty$, but actually we have a decay
of $1/N$.
To show this, note that 
$ f(\bar X_i,\bar X_j) 
 - \int_{\mathbf R} f(\bar X_i, y) u_t(y) dy$
 is mean-zero, so for $j_1\neq j_2$,
\begin{equ}[e:mech]
 \E \Big[
\Big( f(\bar X_i,\bar X_{j_1}) 
 - \int_{\mathbf R} f(\bar X_i, y) u_t(y) dy\Big)
 \Big( f(\bar X_i,\bar X_{j_2}) 
 - \int_{\mathbf R} f(\bar X_i, y) u_t(y) dy\Big)
 \Big] =0
 \end{equ}
 since $\bar X_{j_1}$ and $\bar X_{j_2}$ are independent.

To bound the other terms on the RHS of \eqref{e:sznit},
we need some assumptions on $f$, for instance we assume that $f$
is globally Lipschitz. Then we can
easily bound 
the other terms on the RHS of \eqref{e:sznit} by $|X_j - \bar X_j|$.
Upon taking expectation, and by symmetry of law,
these are just $\E |X_i - \bar X_i|$.
We could then apply Gronwall's lemma to show that
$\E \sup_{t\in [0,T]} |X_i(t) - \bar X_i (t)| \sim 1/\sqrt{N}$.

The key of the above argument is the decay in \eqref{e:1N},
which relies on the mechanism demonstrated in \eqref{e:mech} for 
$L^2$ products of independent and mean-zero random variables.
One also needs some way to ``move'' other terms from the RHS to the LHS
using Gronwall's lemma or some other methods.
We will see that the estimates in singular SPDEs will be much more complicated,
but these basic mechanisms will still be present.
For instance, our Theorem~\ref{th:1} below can be viewed as a result of this flavor, in an SPDE setting,
and in fact the starting point of our proof is indeed close in spirit to  \cite[Sec~I(1)]{MR1108185}
where one subtracts $X_i$ from $Y_i$ to cancel the noise and then bound a suitable norm of the difference.

\section{Da Prato--Debussche and Mourrat--Weber arguments for dynamical $\Phi^4_2$}\label{sec:Phi42}

Before discussing large $N$ problems, we first briefly review
the solution theory for the stochastic quantization of $\Phi^4$ on $\mathbf T^2$, i.e. 
the case $N=1$. 
For simplicity here and in the next section we set $\lambda=1$.
The equation \eqref{e:Phi4-formal} is only formal
in $d=2$ since $\Phi$ is expected to be distributional.
To give a meaning to it, one should regularize and
 insert renormalization terms (also often called ``counter-terms" in the context of
quantum field theory  \footnote{The renormalization constant here for our SPDE is consistent with the one in QFT which is well-known in physics.}). Namely,
we take a sequence of mollified noises $\xi_\eps \to \xi $ (for instance convolving $\xi$ with smooth functions $\chi_\eps$ which approximate Dirac distribution as $\eps\to 0$) and consider the mollified equation
 \begin{equation} \label{e:RPhi4eps}
 \partial_t \Phi_\eps = \Delta \Phi_\eps
 - ( \Phi_\eps^3  - C_\eps \Phi_\eps) + \xi_\eps
\end{equation}
where   $C_\eps$ {\it diverges} as $\eps\to 0$ at a suitable rate,
with some initial condition   $\Phi_\eps(0)=\phi$.
\footnote{We will assume sufficiently regular initial condition $\phi$
in this section and focus on roughness of the equation itself for simplicity.}
If the sequence of constants $C_\eps$ is suitably chosen,
the sequence of smooth solutions $\Phi_\eps$ of \eqref{e:RPhi4eps} will converge to a nontrivial limit
$\Phi = \lim_{\eps\to 0} \Phi_\eps$.

\subsection{Renormalization}
\label{sec:Renormalization}
The following argument is originally by Da Prato--Debussche \cite{DD03}.
Write $\LL = \partial_t - \Delta$ and let 
$Z_\eps$ be the stationary solution 
to  the mollified linear stochastic heat equation $\LL Z_{\eps} = \xi_{\eps}$. The key observation is that the most singular part of $\Phi$ is $Z$, so if we decompose
\begin{equ}[e:DD1]
\Phi_\eps = Z_\eps + Y_\eps,
\end{equ}
we can expect the remainder $Y_\eps$ to converge in a space of better regularity.
Substituting \eqref{e:DD1} into \eqref{e:RPhi4eps} yields
\begin{equ} [e:v-equ]
\LL Y_{\eps}
=
-Y_{\eps}^{3}
-3Z_{\eps}Y_{\eps}^{2}
-3  \cZ_\eps^{\<2>} Y_{\eps}
-\cZ_{\eps}^{\<3>} \;.
\end{equ}
Here, we have defined
\begin{equ}
\cZ_{\eps}^{\<2>}= Z_{\varepsilon}^2-C_\eps  \;,
\qquad
\cZ_{\eps}^{\<3>} =
Z_{\varepsilon}^3- 3C_\varepsilon Z_{\varepsilon}   \;.
\end{equ}

 We will choose $C_\varepsilon=\mathbf{E}[ Z_{\varepsilon}^2(0,0)]$
and the reason will be clear below.
The equation \eqref{e:v-equ} now looks more promising since the noise $\xi_{\eps}$ (whose limit is rough) has dropped out.
This manipulation has not solved the problem of multiplying distributions, since the limit of $Z_\eps$ is still a distribution,
which is why we need renormalization when define its powers.
However $Z_{\eps}$ is a rather concrete object since it is {\it Gaussian} distributed. This makes it possible to study the renormalization and convergence of $Z_{\eps}^{2}$ and $Z_{\eps}^{3} $ via probabilistic methods.
For example, consider the expectation
\begin{equ}[e:Eu2]
\mathbf E [Z_{\eps}(t,x)^{2}]
= \int P(t-s,x-y) P(t-r,x-z) \chi_\eps^{(2)}(s-r,y-z) \,dsdydrdz
\end{equ}
where $\chi_\eps^{(2)}$  is the convolution of the mollifier $\chi_\eps$  with itself and $P$ is the heat kernel.
Due to the singularity of the heat kernel $P$ at the origin, this integral diverges like $O(\log \eps)$ as $\eps\to 0$ in two spatial  dimensions. 
The  constant $C_\eps$
defined above simply subtracts this divergent mean  from $Z_{\eps}^{2}$.
Similarly for the other term $Z_{\eps}^{3}$.
This amounts to considering the renormalized $\Phi^4$ equation \eqref{e:RPhi4eps}.

It can be shown that $Z_\eps$, $\cZ_{\eps}^{\<2>} $ and
$\cZ_{\eps}^{\<3>}$ do converge to nontrivial limits in probability 
as space-time distributions of regularity $\bC^{-\kappa}$ for any $\kappa>0$.
\footnote{Here we denote by $\bC^{\alpha}$ the H\"older--Besov spaces, see e.g. \cite[Appendix~A]{SSZZ2d} for the definitions of these spaces.}
In fact, thanks to Gaussianity of $Z_\eps$, given a smooth test function
$f$ one can explicitly compute any probabilistic moment of
$ \cZ_{\eps}^{\<2>}  (f)$ and prove its convergence.
By choosing $f$ from a  suitable set of wavelets or Fourier basis,
one can apply a version of  Kolmogorov's theorem 
to  prove that $Z_\eps$, $\cZ_{\eps}^{\<2>} $ and
$\cZ_{\eps}^{\<3>}$ converge in  $\bC^{-\kappa}$.
We denote these limits by $Z$, $\cZ^{\<2>} $ and
$\cZ^{\<3>}$. 

\subsection{Fixed point argument}
\label{sec:Fixed point}
Passing \eqref{e:v-equ} to the limit, we get
\begin{equ} [e:v-equ-R]
\LL Y
=
- \Big(Y^{3}+3 Z Y^{2}
+3 \cZ^{\<2>}  Y + \cZ^{\<3>}\Big)\;.
\end{equ}
We can prove local well-posedness of this equation as a {\it classical} PDE, by a standard fixed point argument. For this, we use a classical result in harmonic analysis  under the name ``Young's theorem'' (e.g. \cite[Proposition~4.14]{Hairer14}, \cite[Theorems~2.47 and 2.52]{MR2768550})
which states that if $f \in \mathcal C^{\alpha}$,
$g \in \mathcal C^{\beta}$, and $\alpha+\beta>0$, then $f\cdot g \in \mathcal C^{\min(\alpha, \beta)}$. Thus if we assume that $Y \in \mathcal C^{\beta}$ for $\beta \in (0,2)$, then the worst term in the parenthesis in  \eqref{e:v-equ-R} has regularity $-\kappa$.
By the classical Schauder estimate which states that the heat kernel improves regularity by $2$, the following fixed point map
\begin{equ} [e:fix-pt-Phi42]
Y \mapsto
\LL^{-1}
\Big(Y^{3}+3Z Y^{2}
+3 \cZ^{\<2>} Y + \cZ^{\<3>} \Big)
\end{equ}
is well-defined on $\bC^{\beta}$, since $\bC^{-\kappa+2}\subset \bC^{\beta}$.
With a bit of extra effort, one can show that over short time interval the fixed point map is contractive and thus has a fixed point in $\bC^{\beta}$, which is the solution.
To conclude, one has $\Phi = Z+Y$, which is the local solution to the renormalized $\Phi^4$ equation on $\mathbf T^2$; the above renormalization is well-known as ``Wick renormalization'' and a commonly used notation is to write the equation as
 \begin{equation} \label{e:Phi4wick}
 \partial_t \Phi = \Delta \Phi
 - \Wick{ \Phi^3 } + \xi \;.
\end{equation}

\subsection{A priori estimate}
\label{sec:A priori estimate}
The global (long-time) bound
is first established by Mourrat and Weber \cite{MW17}.
The start point is a {\it PDE energy estimate} -- here it will be an $L^p$ estimate
for $Y$.  
Namely, we multiply \eqref{e:v-equ-R} by $Y^{p-1}$ and integrate over space:
\begin{equ}[e:Lp]
\frac1p  \partial_t  \|Y_t\|_{L^p}^p  + (p-1) \langle \nabla Y, Y^{p-2} \nabla Y \rangle
+ \|Y^{p+2}\|_{L^1} = -\langle Y^{p-1}, 3Z Y^{2}
+3 \cZ^{\<2>} Y + \cZ^{\<3>} \rangle
\end{equ}
where $\langle \;,\;\rangle$ denotes the $L^2$ product.  Note that
the second term on the LHS arises from $\Delta Y$,
and the third term on the LHS arises from $-Y^3$ in  \eqref{e:v-equ-R} (this negative sign is crucial, since it yields a positive quantity on the LHS).

One then needs to control the terms on the RHS. We demonstrate the 
basic idea
with the first term, namely $-\langle Y^{p-1}, 3Z Y^{2}\rangle$.
Standard Besov space inequalities allow us to show that
$$
\Big| 
\langle Y^{p-1}, Z Y^{2}\rangle
\Big|
=\Big| 
\langle Y^{p+1}, Z \rangle
\Big|
\lesssim
\|Y^{p+1}\|_{B_{1,1}^\alpha} \|Z\|_{B_{\infty,\infty}^{-\alpha}}\;.
$$
Recall that $Z \in \bC^{-\kappa}$ (which is essentially the same as
the Besov space $B_{\infty,\infty}^{-\kappa}$),
so $\|Z\|_{B_{\infty,\infty}^{-\alpha}}$ is well-bounded.

One then has the following interpolation inequality
$$
\|Y^{p+1}\|_{B_{1,1}^\alpha}
\lesssim
\|Y^{p+1}\|_{L^1}^{1-\alpha}
\| \nabla (Y^{p+1})\|_{L^1}^\alpha \;.
$$
We would like to use Young's inequality to write the product on the RHS into a sum, and then ``absorb'' them by the 2nd and 3rd terms 
on the LHS of \eqref{e:Lp}. To this end, some manipulation is needed:
$$
\|Y^p \nabla Y\|_{L^1}^\alpha \lesssim 
\|(Y^{\frac{p-2}{2}} \nabla Y )^2\|_{L^1}^{\frac{\alpha}{2}}
\|Y^{p+2} \|_{L^1}^{\frac{\alpha}{2}}
$$
$$
\|Y^{p+1}\|_{L^1}^{1-\alpha} \lesssim 
\|Y^{p+2}\|_{L^1}^{\frac{p+1}{p+2}(1-\alpha)} \;.
$$
Now we see that the terms on  the LHS of \eqref{e:Lp} show up here.
By Young's inequality $ab\le a^{q_1}+b^{q_2}$ with proper choice of $q_1,q_2$ such that $q_1^{-1}+q_2^{-1}=1$,
we can then bound
$\Big| 
\langle Y^{p-1}, Z Y^{2}\rangle
\Big| $  by (say) $\frac{1}{10}$  times the
2nd and 3rd terms 
on the LHS of \eqref{e:Lp} plus a finite constant.
The other terms on the RHS of \eqref{e:Lp} are bounded in the same way.
Putting these bounds together, one obtains global control on the $L^p$ norm.

We emphasize that an energy identity of the form \eqref{e:Lp}
followed by various bounds will be the start-point of our analysis below 
in the large N problems.

\section{Mean field limit of dynamical linear $\sigma$-model in 2D}
\label{sec:MF2D}

In this section, we consider the $N$-vector generalization of \eqref{e:Phi4wick}, which is the stochastic quantization of the linear $\sigma$-model in 2D:
\begin{equation}\label{eq:Phi2d}
\LL \Phi_i=-\frac{1}{N}\sum_{j=1}^N \Wick{\Phi_j^2\Phi_i}+\xi_i,\qquad \Phi_i(0)=\phi_i.
\end{equation}
Note that the solution $\Phi_i$ depends on $N$, but we omit this dependence in our notation.
The solution theory 
by \cite{DD03}  and \cite{MW17} explained in the previous section
corresponds to the case $N=1$.  For a fixed $N$, these well-posedness results are easy to generalize to the present setting.  Our primary goal in this section is to explain a priori bounds which are stable with respect to $N$, and we will then send $N$ to infinity to obtain the following {\it mean field limit}:
\begin{equation}
\LL \Psi_{i}=- \Wick{\E[\Psi_{i}^{2}]\Psi_{i}}+\xi_{i}, \qquad \Psi_i(0)=\psi_i\label{eq:formPsi}.
\end{equation}
The precise meaning of this equation will be given by \eqref{eq:Psi2}
and \eqref{e:1:X} below.

This equation is reminiscent to the mean field limit SDE \eqref{e:barXSDE} 
discussed in the previous section, in the sense that
the equation is self-consistently defined which involves {\it the law of the solution} (here it depends on the second moment of $\Psi_i$). 
On the formal level this equation arises naturally: assuming the initial conditions $\{\phi_{i}\}_{i=1}^{N}$ are exchangeable, 
\footnote{Recall that this means that
their joint probability distribution does not change under permutations of the components.}
the components $\{\Phi_{i}\}_{i=1}^{N}$ will have identical laws, so that replacing the empirical average $\frac{1}{N}\sum_{j=1}^{N}\Phi_{j}^{2}$ by its mean and re-labelling $\Phi$ as $\Psi$ leads us to \eqref{eq:formPsi}. 

 In two space dimensions, \eqref{eq:formPsi} is a {\it singular} SPDE where the ill-defined non-linearity also requires a renormalization. 
 As far as we know this is the first example of singular ``McKean--Vlasov SPDE''
which requires renormalization.
 
The above result is precisely given in the next theorem.

\begin{theorem}[Large $N$ limit of the dynamics for $d=2$]\label{th:1}
Let $\{(\phi_{i}^{N},\psi_{i}) \}_{i=1}^{N}$ be  random initial datum in $\bC^{-\kappa}$ for some small $\kappa>0$ and all moments finite.  Assume that for each $i \in \N$,  $\phi_{i}^{N}$ converges to $\psi_{i}$ in $L^{p}(\Omega; \bC^{-\kappa})$ for all $p>1$, $\frac1N\sum_{i=1}^N\|\phi_i^N-\psi_i\|_{\bC^{-\kappa}}^p$ converges in probability to $0$ and $(\psi_{i})_{i}$ are iid. 

Then for each component $i$ and all $T>0$, the solution $\Phi_{i}$ to \eqref{eq:Phi2d} converges in probability to $\Psi_i$ in $C([0,T], \bC^{-1}(\mathbb T^2))$ as $N\to \infty$, where $\Psi_i$ is the unique solution to the mean-field SPDE formally described by
\begin{equation}\label{eq:Psi2}
\LL \Psi_i= -\mathbf{E}[\Psi_i^2-\tilde Z_i^2] \Psi_i+\xi_i,\quad \Psi_i(0)=\psi_i,
\end{equation}
and $\tilde Z_i$ is the stationary solution to $\LL \tilde Z_i=\xi_i$.
Furthermore, if $(\phi_{i}^{N},\psi_{i})_{i=1}^{N}$ are exchangeable, then for each $t>0$ it holds that
\begin{equation}
\lim_{N \to \infty}\E\|\Phi_i(t)-\Psi_i(t)\|_{L^2(\mathbb{T}^2)}^2=0.
\end{equation}
\end{theorem}

Note that \eqref{eq:Psi2}
is still ``formal'' i.e. its meaning needs interpretation, see \eqref{e:1:X} below.
We could actually prove this convergence result under more general conditions for initial data (see \cite[Sec.~4]{SSZZ2d}).

\begin{remark}\label{rmk:1d}
In one space dimension, the model does not need any renormalization,
so the result is much simpler to understand.
Consider the  equation on $\R^+\times \mathbb{T}$:
\begin{equation}\label{eq:1}\LL \Phi_i=-\frac{1}{N}\sum_{j=1}^N\Phi_j^2\Phi_i+\xi_i, \quad \Phi_i(0)=\phi_i,
\end{equation}
with $\LL=\p_t-\Delta+m$ for $m\geq0$, $1\leq i\leq N$.
We assume $\mathbf{E}\|\phi_i\|_{L^2}^2\lesssim1$ uniformly in $i, N$.
In $d=1$, we can show that
 the limiting equation as $N\to \infty$ is given by
\begin{equs}\label{eq:Psi}
\LL\Psi_i =-\mu\Psi_i+\xi_i,
\qquad
\mbox{with} \qquad \mu(t,x)=\mathbf{E}[\Psi_i(t,x)^2] \;.
\end{equs}
Namely
 $\frac{1}{N}\sum_{j=1}^N\Phi_j^2$ in the equation  self-average in the large $N$ limit to its mean.
See  \cite[Supplementary Material]{SSZZ2d}.
\end{remark}

\subsection{Uniform in $N$ bounds}
Following the  trick of Da-Prato and Debussche discussed in Sec.~\ref{sec:Phi42}, we consider the decomposition $\Phi_i = Z_i + Y_i$, where $Z_{i}$ is a solution to the linear SPDE
\begin{equation}\label{eq:li1}
\LL Z_i=\xi_i,\quad Z_i(0)=z_{i},
\end{equation}
and $Y_{i}$ solves the remainder equation
\begin{equation}\label{eq:22}
\LL Y_i =-\frac{1}{N}\sum_{j=1}^N(Y_j^2Y_i+Y_j^2Z_i+2Y_jY_iZ_j+2Y_j\Wick{Z_iZ_j}
+ \Wick{Z_j^2}Y_i+\Wick{Z_i Z_j^2 })
\end{equation}
with some initial condition $Y_i(0) =y_i$.
Here we assume that the initial datum satisfy  $\mathbf{E}\|z_i\|_{\bC^{-\kappa}}^p\lesssim1$ for $\kappa>0$ small enough and every $p>1$, and $\mathbf{E}\|y_i\|_{L^2}^2\lesssim1$, uniformly in $i, N$.
 \eqref{eq:22} has similar structure as
\eqref{e:v-equ}: it has on the RHS terms cubic in $Y$, in $Z$, and cross-terms between $Y$ and $Z$. 
 
The notation $\Wick{Z_iZ_j} $, $ \Wick{Z_j^2}$ and $\Wick{Z_i Z_j^2}$ denotes the Wick products. It is easier to explain the Wick products for 
$\tilde Z_i$ which is the stationary solution to $\LL \tilde Z_i=\xi_i$.
As in Section~\ref{sec:Renormalization},
and using independence of $(\tilde Z_i)_{i=1}^N$,
we should define

\begin{equ}[e:wick-tilde]
\Wick{\tilde{Z}_i\tilde{Z}_j}=
\begin{cases}
\lim\limits_{\varepsilon\to0}(\tilde{Z}_{i,\varepsilon}^2-a_\varepsilon)  &  (i=j)\\
 \lim\limits_{\varepsilon\to0}\tilde{Z}_{i,\varepsilon}\tilde{Z}_{j,\varepsilon} & (i\neq j)
\end{cases}
\quad
\Wick{\tilde{Z}_i\tilde{Z}_j^2} =
\begin{cases}
 \lim\limits_{\varepsilon\to0}(\tilde{Z}_{i,\varepsilon}^3-3a_\varepsilon \tilde{Z}_{i,\varepsilon})   & (i=j)\\
 \lim\limits_{\varepsilon\to0}(\tilde{Z}_{i,\varepsilon}\tilde{Z}_{j,\varepsilon}^2-a_\varepsilon \tilde{Z}_{i,\varepsilon}) & (i\neq j)
\end{cases}
\end{equ}
 where $a_\varepsilon=\mathbf{E}[\tilde Z_{i,\varepsilon}^2(0,0)]$ is a divergent constant independent of $i$
 (in fact it is equal to $C_\eps$ in  Section~\ref{sec:Renormalization}).
 As in the previous section, the above limits exist in 
 $\bC^{-\kappa}$
  for $\kappa>0$.
 The Wick products $\Wick{Z_iZ_j} $, $ \Wick{Z_j^2}$ and $\Wick{Z_i Z_j^2}$
are just suitable modifications by additional terms involving the initial condition $z$ of $Z$ (also subtracting the same $a_\varepsilon$, see \cite[Sec.~2.1]{SSZZ2d} for details).

A key step in proof of
Theorem \ref{th:1} is a new uniform in $N$ bounds through suitable energy estimates on the remainder equation \eqref{eq:22}.  This is inspired by the approach of a priori bounds in \cite{MW17} discussed in Sec.~\ref{sec:A priori estimate}, but  subtleties arise as we track carefully the dependence of the bounds on $N$. 
 Indeed, the   approach  to obtain global in time bounds for  $N=1$ is to exploit the damping effect from
the term $-Y^{3}$  in \eqref{e:v-equ-R}, which corresponds to the term 
  $-Y_j^2Y_i$.
However, the extra factor $1/N$ makes this effect weaker as $N$ becomes large.  In fact, the moral is that we cannot exploit the strong damping effect at the level of a fixed component $Y_{i}$, rather we're forced to consider aggregate quantities, and ultimately we focus on the empirical average of the $L^2$-norm (squared) instead of the $L^p$-norm for arbitrary $p$ as in Sec.~\ref{sec:Phi42}. 
We skip the details in this expository note and jump into the discussion on the mean-field equation and the proof of convergence.

\subsection{The mean-field SPDE}
We now discuss a bit more the solution theory for the mean-field SPDE \eqref{eq:Psi2}.   While the notion of solution we use is again via the Da-Prato/Debussche trick, the well-posedness theory requires more care than that for $\Phi^4_2$  discussed in Sec.~\ref{sec:Phi42}, since here the equation depends on its own law and we cannot proceed by pathwise arguments alone. 
Again, we understand \eqref{eq:Psi2} via the decomposition $\Psi_i=Z_i+X_i$ with $X_i$ satisfying
 \begin{equ}[e:1:X]
\LL X_i =-(\E[X_j^2]X_i+\E [X_j^2]Z_i+2\E[X_jZ_j]X_i+2\E [X_jZ_j]Z_i
) + (\cdots).
\end{equ}
This follows by a bit algebra, noting that $Z_i^2$ and the $\tilde Z_i^2$
in  \eqref{eq:Psi2} ``basically'' cancel,
except that the latter one is stationary while the former one is not 
-- this is the term $(\cdots)$ which we omit in the following discussion. 
What's more important is that here we actually introduced an independent copy $(X_j,Z_j)$ of $(X_i,Z_i)$ inside $\E$, which turns out to be useful for both the local and global well-posedness of \eqref{eq:Psi2}.  Indeed, one point is that the term $\E [X_jZ_j]Z_i$ in \eqref{e:1:X} cannot be understood in a classical sense; however  we can view it as a conditional expectation $\E [X_jZ_jZ_i|Z_i]$ and use properties of the Wick product $Z_{i}Z_{j}$ (like in \eqref{e:wick-tilde}) to give a meaning to this. 
  In fact, after taking expectation, $\E[X_j^2]X_i$ in \eqref{e:1:X} also plays the role of the damping mechanism, which helps us to obtain uniform bounds on the mean-squared $L^2$-norm of $X_i$.
  
 Assume the corresponding decomposition for initial datum
 $\psi_i=z_i+\eta_i$. We denote 
 $L^{2}_TL^{2} =L^2 ([0,T], L^2)$.
 We have:

\bl\label{lem:l2} There exists a universal constant $C$ such that

\begin{align}
\sup_{t\in[0,T]}\E \|X_{i}\|_{L^{2}}^{2}+\E \|\nabla X_{i}\|_{L^{2}_TL^{2} }^{2}+\| \E X_{i}^{2} \|_{L^{2}_T L^{2} }^{2}+m\E \| X_{i}\|_{L^{2}_T L^{2}   }^{2}
 \leq C \int_{0}^{T}R \dif t+\E\|\eta_i\|_{L^2}^2. \label{s30}
\end{align}
\el
Here, $R$ contains terms which only depend on renormalized powers of $Z$ and is well bounded.
It has the following explicit form where $i \neq j$ and $s>0$ is a sufficiently small parameter
\begin{equ}
R\eqdef 1 +\big ( \E\|Z_i\|_{\bC^{-s}}^{2}  \big )^{\frac{2}{1-s}} +\E\|\Wick{Z_{j}^{2}Z_{i}}\|_{\bC^{-s}}^{2}
 +C\big ( \E \|\Wick{Z_{j}Z_{i}}\|_{\bC^{-s}}^{2} \big )^{2} +C\big (\E \|\Wick{Z_{i}^{2}} \|_{\bC^{-s}} \big )^{4}\;.
\end{equ}

Recall in Sec.~\ref{sec:Phi42} that the start point of the global estimate is a energy identity \eqref{e:Lp}.
The proof of Lemma~\ref{lem:l2}
is,  again similarly as in Sec.~\ref{sec:Phi42}, based on an energy identity
	\begin{align}
	\frac{1}{2}\frac{\dif}{\dif t}\E \|X_{i}\|_{L^{2}}^{2}+\E \|\nabla X_{i}\|_{L^{2}}^{2}+ \|\E X_{i}^{2}\|_{L^{2}}^{2}+m\E\|X_i\|_{L^2}^2 =I^{1}+I^{2}+I^{3}, \label{s20}
	\end{align}
where
	\begin{align}\label{eq:I}
	I^{1}\eqdef \E\langle X_{i},\Wick{ Z_{i}Z_{j}^{2}} \rangle, \qquad
	I^{2}\eqdef \E\langle X_{i}^{2},  \Wick{Z_{j}^{2}}\rangle+2\E\langle X_{i}X_{j}, Z_{i}Z_{j}\rangle, \qquad
	I^{3}\eqdef 3\E\langle X_{i}^{2}X_{j},Z_{j} \rangle .
	\end{align}
We then proceed by bounding $I^{1,2,3}$
by
$\frac{1}{10} \big (\E \|\nabla X_{i}\|_{L^{2}}^{2}+ \|\E X_{i}^{2}\|_{L^{2}}^{2} \big )$
plus terms that only depend on $Z$ (which eventually constitute $R$).
We refer to 
 \cite[Sec.~3.2]{SSZZ2d} for these details.

\subsection{Proof of convergence}
We now explain the basic ideas in the proof of Theorem~\ref{th:1}.  
Recall that  $\Phi_{i}=Y_{i}+Z_{i}$, and  $\Psi_{i}=X_{i}+Z_{i}$,  we define (in the same spirit as we did in Sec.~\ref{sec:SDE})
\begin{equation}
v_{i} \eqdef Y_{i}-X_{i} \nonumber.
\end{equation}
Now we have the following energy identity
	\begin{align}
	& \frac12\frac{\dif}{\dif t}\sum_{i=1}^{N} \|v_{i}\|_{L^{2}}^{2}+ \sum_{i=1}^{N}\|\nabla v_{i}\|_{L^{2}}^{2}+m\sum_{i=1}^{N} \|v_{i}\|_{L^{2}}^{2}+\frac{1}{N}\sum_{i,j=1}^{N} \|Y_{j}v_{i}\|_{L^{2}}^{2}+\frac{1}{N}\bigg \|\sum_{j=1}^{N}X_{j}v_{j}  \bigg \|_{L^{2}}^{2} = \sum_{k=1}^3 I_k^N \label{diff11}
	\end{align}
	where
	\begin{align}
	I_{1}^{N} & \eqdef
	-\frac{1}{N}\sum_{i,j=1}^{N}  \big ( 2\langle v_{i}v_{j},\Wick{ Z_{j}Z_{i}}\rangle+\langle v_{i}^{2},\Wick{Z_{j}^{2}}\>+2\langle v_{i}^{2}Y_{j},Z_{j}\rangle \big )\;,
	\nonumber \\
	I_{2}^{N}& \eqdef
	-\frac{1}{N}\sum_{i,j=1}^{N} \langle v_{i}v_{j}, \big (X_{i}Y_{j}+(3X_{j}+Y_{j} )Z_{i}  \big )\rangle\;,
	\nonumber \\
	I_{3}^{N}&\eqdef
	-\frac{1}{N}\sum_{i,j=1}^N \Big\langle \big [\Wick{Z_{j}^{2} }-\E\Wick{Z_{j}^{2} }+ X_{j}(X_{j}+2Z_{j} )-\E  X_{j}(X_{j}+2Z_{j})  \big ] (X_{i}+Z_{i})\;,\;v_{i} \Big \rangle\;,
	\end{align}
	In the definition of $I_{3}^{N}$, to have a compact formula, we slightly abuse notation for the contribution of the diagonal part $i=j$, where we understand $Z_{i}Z_{j}$ to be $\Wick{Z_{i}^{2}}$ and $\Wick{Z_{j}^{2}}Z_{i}$ to be $\Wick{Z_{i}^{3}}$. 
The above could be viewed as a much more complicated version of \eqref{e:Lp}.
	
Next, we bound $I_{1,2,3}^{N}$ one by one.
This requires substantial calculations,
but we illustrate a key point with $I_{3}^{N}$ in the following.
Define
	\begin{equation}
	G_{j}\eqdef  (X_{j}^2-\E X_j^2)+2(X_jZ_{j}-\E X_{j}Z_{j})+(\Wick{Z_{j}^{2} }-\E\Wick{Z_{j}^{2} })\eqdef G_j^{(1)}+G_j^{(2)}+G_j^{(3)}.\nonumber
	\end{equation}
It is not very hard (see \cite[(4.18)]{SSZZ2d}) to prove the following bound, for $s>0$ small:
	\begin{align}\label{eq:zz}
	I_3^N & \leq C(\bar{R}_N+\bar{R}_N^\prime)+\frac18\frac1N \Big\|\sum_{i=1}^NX_iv_i\Big\|_{L^2}^2+\frac{1}{4}\sum_{i=1}^N\| v_i\|_{H^1}^2
	\\&+C\Big(\sum_{i=1}^N\|v_i\|_{L^2}^2\Big)
	\Big[R_N^Z+1
	+\frac{1}{N}\sum_{i=1}^N\Big(\|X_i\|_{L^4}^{4/(1-2s)}+\|\Lambda^sX_i\|_{L^4}^4\Big)\Big],\nonumber
	\end{align}
	with $\Lambda=(1-\Delta)^{\frac12}$ and
	\begin{align*}
	\bar{R}_N & \eqdef \frac{1}{N}\Big\|\sum_{j=1}^NG_j^{(1)}\Big\|_{H^{s}}^2
	+\!\!\! \sum_{k\in\{2,3\}}\frac{1}{N}\Big\|\sum_{j=1}^N G_j^{(k)}\Big\|_{H^{-s}}^2 ,
	\qquad
	\bar{R}_N^\prime  \eqdef \sum_{k\in\{2,3\}} \frac{1}{N^2}\sum_{i=1}^N \| \sum_{j=1}^N G_j^{(k)} Z_i\|_{H^{-s}}^2,
	\\
	R_N^Z & \eqdef \Big(\frac{1}{N}\sum_{i=1}^N\|\Lambda^{-s}Z_i\|_{L^\infty}^2\Big)^{\frac{1}{1-s}}.
	\end{align*}
A key point is that $G_j^{(k)}$ are centered. 
Also, components of $X$ are independent, 
and so are components of $Z$.
Therefore, although 
the terms in $\bar{R}_N$, $\bar{R}_N^\prime$
appear to be ``order $O(N)$'',
we can actually show that upon taking expectation
they are bounded uniformly in $N$.
Essentially,
this boils down to 
the following fact:
for mean-zero independent random variables $U_{1},\dots,U_{N}$ taking values in a Hilbert space $H$, we have (``gaining a factor $1/N$'')
\begin{equation}
\E \Big \| \sum_{i=1}^{N}U_{i} \Big \|_{H}^{2} =\E \sum_{i=1}^{N}\|U_{i}\|_{H}^{2}.\label{eq:Ui}
\end{equation}
This is ultimately similar with the mechanism \eqref{e:mech}
in the earlier section for the ODE case.
After considerably technical estimates, we can show that 
$\sup_{t\in[0,T]}\|v_i(t)\|_{L^2}^2\to0$ in probability, as $N\to \infty$,
and Theorem~\ref{th:1} holds 
(see \cite[Sec.~4]{SSZZ2d} for details).

\section{Convergence of (Wick renormalized) linear $\sigma$-model to Gaussian free field in 2D}\label{sec:GFF2d}

In this section, we are only interested in the stationary setting,
and we simply write
$Z$ for the  stationary solution to the  linear equation $\LL Z=\xi$.

We now study the invariant measure for the mean field equation \eqref{eq:Psi2}, namely
\begin{equation}\label{eq2:Psi}
\LL\Psi=-\mathbf{E}[\Psi^2- Z^2]\Psi+\xi,  
\end{equation}
with $\mathbf{E}[\Psi^2-Z^2]=\mathbf{E}[X^2]+2\mathbf{E}[XZ]$ for $X=\Psi-Z$.
The remainder
$X$ satisfies
\begin{equation}\label{eq2:X}
\LL X=-\mathbf{E}[X^2+2XZ](X+Z),\quad X(0)=\Psi(0)-Z(0).
\end{equation}

An immediate observation is that $Z$ is a stationary solution to  \eqref{eq2:Psi}.  
This follows since the unique solution to \eqref{eq2:X} starting from zero is identically zero (or one could simply set $\Psi=Z$ in \eqref{eq2:Psi}).  
Since the (massive) Gaussian free field  $\mathcal N(0,\frac12(-\Delta+m)^{-1})$ is invariant measure for the $Z$ equation,
we see that Gaussian free field is an invariant measure for  \eqref{eq2:Psi}.

\subsection{Uniqueness of invariant measure for mean-field equation}
Below, we use a semigroup $P_t^*\nu$ to denote the law of $\Psi(t)$ with the initial condition distributed according to a measure $\nu$.   A natural question is whether \eqref{eq2:Psi} has unique  invariant measure.
Since the equation involves its own law, 
it is unclear if the general ergodic theory can be applied directly in this setting. Fortunately, it has a strong damping property in the mean-square sense, which comes to our rescue and allows us to proceed directly by a priori estimates.

\begin{remark}
If we assume that $\mathbf{E}[\Psi^2- Z^2] =:\mu$ is simply a constant,
then uniqueness is a simple calculation.
Note that the limiting equation $\LL\Psi=-\mu\Psi+\xi$ has a Gaussian invariant measure $\mathcal N(0,\frac12(-\Delta+m+\mu)^{-1})$.
	Assuming $\Psi \sim \mathcal N(0,\frac12(-\Delta+m+\mu)^{-1})$ and $Z  \sim \mathcal N(0,\frac12(-\Delta+m)^{-1})$,
	the self-consistent condition $\mathbf{E}[\Psi^2- Z^2] = \mu$ then yields
	\begin{equ}
	\frac12\sum_{k\in \mathbb Z^2} \Big(\frac{1}{|k|^2+m+\mu}-\frac{1}{|k|^2+m}\Big) = \mu\;.
	\end{equ}
For $\mu+m\ge 0$ we only have one solution $\mu=0$, since the LHS is monotonically decreasing in $\mu$.

If $d=1$ as in Remark~\ref{rmk:1d} and $\mu$ is constant, then it is the unique solution to $\frac12\sum_{k\in \mathbb Z}
\frac{1}{k^2+m+\mu}= \mu$.
\end{remark}

More precisely, from the basic energy identity \eqref{s20},
we can prove that (\cite[Lemma~5.3]{SSZZ2d})
	\begin{equation}
	\frac{1}{2}\frac{\dif}{\dif t} \E \| X_{i} \|_{L^2}^2+\frac{1}{2}\E\|\nabla X_i\|_{L^2}^2 +m\E \|X_{i}\|_{L^2}^2+ \|\E X_{i}^2\|_{L^2}^2 \lesssim 1.\label{se52}
	\end{equation}
	As a consequence, noting that $ \|\E X_{i}^2\|_{L^2}^2 \geq \big ( \E \| X_{i} \|_{L^2}^2 \big )^{2}$, it is not hard to deduce
	\begin{align}
	\sup_{t>0}(t\wedge1)\mathbf{E}[\|X_{i}(t)\|_{L^2}^2] \lesssim1,\label{bd2:uniX}
	\end{align}
	where the implicit constant is independent of the initial data.  
	On the other hand, from  \eqref{s20} one can also prove
	\begin{equation}
	\frac{\dif}{\dif t}\E\|X_{i}\|_{L^{2}}^{2}+m\E\|X_{i}\|_{L^{2}}^{2} \leq C_0(\E \|\Wick{Z_{2}Z_{1}}\|_{\bC^{-s}}^{2}+(\E\|Z_1\|_{\bC^{-s}}^2)^{\frac{1}{1-s}}+1)\mathbf{E}\|X_{i}\|_{L^2}^2  =:  m_0 \mathbf{E}\|X_{i}\|_{L^2}^2. \label{se50}
	\end{equation}
	Applying Gronwall's inequality over $[1,t]$ leads to
	\begin{equation}
	e^{(m-\frac{m_0}{2})t}\E \|X_{i}(t)\|_{L^{2}}^{2} \lesssim \E \|X_{i}(1)\|_{L^{2}}^{2} \nonumber,
	\end{equation}
	so if $m$ is sufficiently large, i.e. $m \geq m_0$, using \eqref{bd2:uniX}, and taking the supremum over $t \geq 1$, we see that $\Psi$ approaches $Z$ in $L^2$ sense exponentially fast
	\begin{equation}
\sup_{t \geq 1}e^{\frac{mt}{2} }\mathbf{E}\|\Psi(t)-Z(t) \|_{L^2}^2 \leq C.\label{se51}
\end{equation}
From this bound we can easily deduce that 
 for sufficiently large mass, the unique invariant measure to \eqref{eq2:Psi} is Gaussian.  To this end, define the $\bC^{-1}$-Wasserstein distance
\begin{equ}[e:Wasserstein]
\mathbb{W}_p'(\nu_1,\nu_2):=\inf_{\pi\in\mathscr{C}(\nu_1,\nu_2)}\left(\int\|\phi-\psi\|_{\bC^{-1}}^p\pi(\dif \phi,\dif \psi)\right)^{1/p}
\end{equ}
where $\mathscr{C}(\nu_1,\nu_2)$ denotes the collection of all couplings of $\nu_1, \nu_2$ satisfying $\int\|\phi\|_{\bC^{-1}}^p\nu_i(\dif \phi)<\infty$.

\bt\label{th:u} For sufficiently large $m$  the unique invariant measure to \eqref{eq2:Psi} supported on $\bC^{-\kappa}$ is $\mathcal N(0,\frac12(-\Delta+m)^{-1})$, the law of the Gaussian free field.
\et

\subsection{Large N limit of the linear $\sigma$-model in 2D}

The solutions $(\Phi_i)_{1\leq i\leq N}$ to \eqref{eq:Phi2d} form a Markov process on $(\bC^{-\kappa})^{N}$ which admits a unique invariant measure, henceforth denoted by $\nu^{N}$. 
Denote by $\nu^{N,i}$ the marginal law of the $i$-th component, 
and $\nu^N_k$ the marginal law of the first $k$ components. 
 Our goal in this section is to study the large $N$ behavior of $\nu^{N}$ and show that for sufficiently large mass, as $N \to \infty$, it's marginals are simply products of the Gaussian invariant measure for $\Psi$ identified in Theorem \ref{th:u}.  For this we rely heavily on the estimates  for the remainder $Y$.  It will be convenient to have a stationary coupling of the linear dynamic
 $\LL Z_i=\xi_i$ 
 and the nonlinear dynamic \eqref{eq:Phi2d} respectively, namely,
 there exists a  stationary process $(\Phi_i^N, Z_i)_{1\leq i\leq N}$ such that the components $\Phi_i^N, Z_i$ are stationary  solutions to \eqref{eq:Phi2d} and $\LL Z_i=\xi_i$, respectively.

Write $\mathbb{W}_2$ for the 
 $\bC^{-\kappa}$-Wasserstein distance
which is defined as in 
\eqref{e:Wasserstein} with  $\bC^{-1}$ replaced by  $\bC^{-\kappa}$
and $p$ replaced by $2$.

\bt\label{th:m1} Let $\nu= \mathcal{N}(0,\frac12(m-\Delta)^{-1})$ be the Gaussian free field with covariance $\frac12(m-\Delta)^{-1}$.
For sufficiently large $m$ 
one has
\begin{equation}
\mathbb{W}_2(\nu^{N,i},\nu) \leq CN^{-\frac{1}{2}}. \label{West}
\end{equation}
Furthermore, $\nu^N_k$ converges to $\nu\times...\times \nu$, as $N\to \infty$.
\et

 The stationarity of the joint law of $(\Phi_i, Z_i)$ implies that also $Y_i=\Phi_i-Z_i$ is stationary. 
The key step to prove Theorem~\ref{th:m1} is to show
	\begin{equation}
	\E \|Y_{i}^{N}\|_{H^{1}}^{2} \leq CN^{-1} \label{se8},
	\end{equation}
	which implies \eqref{West} by definition of the Wasserstein metric and the embedding $H^{1} \hookrightarrow \bC^{-\kappa}$ in $d=2$.  By symmetry of law, 
	$$
		\E \|Y_{i}^{N}\|_{H^{1}}^{2}
=\frac1N\sum_{j=1}^N \E\|\nabla Y_j\|_{L^2}^2+\frac1N\sum_{j=1}^N\E\|Y_j\|_{L^2}^2
	$$
Now again using the energy identity for $Y_i$, combined with the stationarity of $(Y_{j})_{j}$ and $(Z_{j})_{j}$, we find
	\begin{align*}
	&\frac1N \sum_{j=1}^N\E\|\nabla Y_j\|_{L^2}^2
	+ \frac{m}{N} \sum_{j=1}^N\E\|Y_j\|_{L^2}^2
	+\frac{1}{N^2}\E\bigg\|\sum_{i=1}^NY_i^2\bigg\|_{L^2}^2
	\\
	& \leq \frac{C}{N}\E R_N^0
	+\E\bigg(\frac1N\sum_{j=1}^N\|Y_j \|_{L^2}^2 \cdot D_N \bigg),
	\end{align*}
for certainly quantities $R_N^0$ and $D_N$ which only depend on  $Z_j$,
$\Wick{Z_iZ_j} $, $ \Wick{Z_j^2}$ and $\Wick{Z_i Z_j^2}$.
For instance,
\begin{equ}\label{eq:R0}
R_N^0=\frac{1}{N^2} \sum_{i=1}^{N}\Big\|\sum_{j=1}^{N} \nabla^{-s}(\Wick{Z_{j}^{2}Z_{i}})\Big\|^2_{L^2}. 
\end{equ}
To show the decay \eqref{se8}, 
we note that $R_N^0$ apparently behaves as $O(N)$ as $N$ becomes large (three sums versus $N^{-2}$). However, it turns out that 
$\E R_N^0 \sim O(1)$. This is because 
$\nabla^{-s}(\Wick{Z_{j}^{2}Z_{i}})$ is centered (mean-zero),
and they are independent for different values of $j$ -- a phenomena already hinted in \eqref{e:1N}.

Unfortunately, $D_N$ needs to be treated differently, since it is not centered.
Setting $A\eqdef\E D_N$, we may re-center $D_{N}$  above:
	\begin{equs}[e:centerD]
	&\frac1N\sum_{j=1}^N\E\|\nabla Y_j(t)\|_{L^2}^2
	+(m-A)\frac1N\sum_{j=1}^N\E\|Y_j(t)\|_{L^2}^2+\frac{1}{N^2}\E\bigg\|\sum_{i=1}^NY_i^2(t)\bigg\|_{L^2}^2
	\\&\leq C\frac1N\E R_N^0+\frac{1}{2}\frac1{N^2}\E\bigg(\sum_{j=1}^N\|Y_j \|_{L^2}^2\bigg)^2+\E|D_N-A|^2.
	\end{equs}
The explicit form of 	$D_N$ is rather complicated so we do not give it here,
but just as what we did for $\E R_N^0$, we can use
centering and 
 independence to show that
\begin{equation}
	\E|D_N(t)-A|^2 \leq \frac{C}{N}. \label{se6}
\end{equation}	
These bounds then yield
\eqref{se8}, provided that $m\geq A+1$ so that the 2nd term
on the LHS is positive.
This is where we need to assume $m$ large in this proof.
(Remark that if we didn't set $\lambda=1$,
then assuming $\lambda>0$ small rather than $m$ large 
our argument will also work.)

\section{Exact correlation formulae for some $O(N)$ invariant observables}
\label{sec:observables}

Now we discuss  {\it observables} in the stationary setting.
In QFT models with continuous symmetries, physically interesting quantities involve more than just a component of the field itself but also
quantities composed by the fields which preserve the symmetries, called invariant observables. For linear $\sigma$ model
a natural quantity  that is invariant under $O(N)$-rotation is the ``length'' of $\Phi$;
another being the quartic interaction in \eqref{e:Phi_i-measure}.
We thus consider the following two $O(N)$ invariant observables, where
$\Phi $ is distributed as $\nu^N$
\begin{equ}[e:twoObs]
\frac{1}{\sqrt{N}} \sum_{i=1}^N\Wick{\Phi^2_i},\qquad 
 \frac{1}{N} \Wick{\Big(\sum_{i=1}^N\Phi_i^2\Big)^2}.
\end{equ}
We establish the large $N$ tightness of these observables as {\it random fields} in suitable Besov spaces by using  iteration to derive  improved uniform estimates in the stationary case.
Note that \cite{MR578040} considered correlations of these observables.
Our SPDE approach allows us to study these observables as {\it random fields} with precise regularity as $N\to \infty$ which is new.

\begin{theorem}\label{th:1.3}
For $m$ large enough, the following result holds for any $\kappa>0$:
\begin{itemize}[leftmargin=.4in]
  \item  $\frac{1}{\sqrt{N}} \sum_{i=1}^N\Wick{\Phi^2_i}$ is tight in $B^{-2\kappa}_{2,2}$, and  $\frac{1}{N} \Wick{(\sum_{i=1}^N\Phi_i^2)^2}$ is tight in $B^{-3\kappa}_{1,1}$.
  \item  
 The Fourier transform of the two point correlation function of $\frac{1}{\sqrt{N}} \sum_{i=1}^N\Wick{\Phi^2_i}$
  in the limit as $N\to \infty$ is
  given by the explicit formula $2\widehat{C^2} /(1+2\widehat{C^2})$, where $C=\frac12(m-\Delta)^{-1}$ and $\widehat{C}$ is the Fourier transform.
 Moreover
 \begin{equ}[e:EPhi4]
  \lim_{N\to \infty} \E \frac{1}{N} \Wick{\big(\sum_{i=1}^N\Phi_i^2\big)^2}
  =
  -4 \sum_{k\in \mathbb Z^2} \widehat{C^2}(k)^2 /(1+2\widehat{C^2}(k)).
 \end{equ}
\end{itemize}
\end{theorem}
The exact expectation formula \eqref{e:EPhi4} above was also derived in \cite[Theorem~3]{MR578040}. In our theorem above
we also showed that the limiting observables can be viewed as random fields in suitable Besov spaces.

This shows that although (for large enough $m$)
the invariant measure of $\Phi_i$ converges
as $N\to \infty$ to the Gaussian free field,
the limits of the observables \eqref{e:twoObs} have different laws than those if $\Phi_i$ in \eqref{e:twoObs}  were replaced by
$Z_i$:
\begin{equ}[e:two-OZ]
	\frac{1}{\sqrt{N}} \sum_{i=1}^N\Wick{Z^2_i},
	\qquad  \frac{1}{N} \Wick{\Big(\sum_{i=1}^NZ_i^2\Big)^2}.
\end{equ}
A simple application of Wick theorem shows that the two-point correlation of the first observable is $2C^2$ and the expectation of the second one is $0$.
We remark that 
expanding 
 $1 /(1+2\widehat{C^2})$
into geometric series,
we indeed get the sum of ``bubble diagrams'' in \eqref{e:bubbles},
as predicted by perturbation theory in physics.

The tightness part of Theorem~\ref{th:1.3} requires establishing upper bounds on the Besov norms of the observables. However, {\it upper bounds} obtained from PDE estimates are clearly not enough to prove {\it exact} correlation formulas in the second part of the theorem. We now focus on the second part of the theorem since it is more interesting, and we only consider the first observable below and we denote its two point correlation function by $G_N$.

\subsection{Dyson--Schwinger equations}
First we need an algebraic step, which proves the following identity
\begin{equ}[e:alg-id]
\quad (1+\frac{2(N+2)}{N} \widehat{C^2}) \widehat{G_N} 
= 2 \widehat{C^2} + \widehat{Q_N}/N\quad
\end{equ}
where
\vspace{-2ex}
\begin{equ}[e:QN]
Q_N (x-z) = -\int C(x-y)C(x-z) \E \Big(\Wick{\Phi_i \sum_j \Phi_j^2  (y)}\Phi_i(z)\Big)dy
 + \mbox{6 pt correlation}
\end{equ}
where we have omitted the explicit form of the extra terms but they are 
essentially some 6-point correlations of $\Phi$.
The proof of \eqref{e:alg-id} relies on iteratively application of Dyson-Schwinger equation.
Write $\PPhi^2 \eqdef \sum_{i=1}^N \Phi_i^2$.
Dyson-Schwinger equation is a hierarchy of equations,
which relate correlation functions of different orders.
For instance, we have
 \begin{equs}
\int & C(x-z) \E\Big[ \Phi_1(x) \Wick{\Phi_1\PPhi^2(z)} \Wick{\PPhi^2(y)} \Big]d z
\\
& =
- \E[ \Wick{\PPhi^2(x)} \Wick{\PPhi^2(y)} ]
+ 2NC (x-y) \E [\Phi_1 (x) \Phi_1(y)]
\end{equs}
which relate 2-point,  4-point,  and 6-point correlations. Also,
\begin{equs}
\frac{1}{N}& \int C(x-z_1)C(x-z_2) \E  \Big[ \Wick{\Phi_1\PPhi^2(z_1)}\Wick{\Phi_1\PPhi^2(z_2)} \Wick{\PPhi^2(y)} \Big]d z_1d z_2
\\
&=-2\int  C(x-y)C(x-z)
 \E \Big[\Wick{\Phi_1 \PPhi^2(z)} \Phi_1(y) \Big]\,d z
\\
&+\frac{N+2}{N} \int C(x-z)^2 \E\Big[ \Wick{\PPhi^2(y)} \Wick{\PPhi^2(z)}  \Big]d z
\\
&\qquad - \int C(x-z) \E\Big[ \Phi_1(x) \Wick{\Phi_1\PPhi^2(z)} \Wick{\PPhi^2(y)} \Big]d z
\end{equs}
which relate 4-point,  6-point,  and 8-point correlations.
Since $\Phi$ is distributional, the rigorous meaning of these equations
should be given via suitable approximation. We refer to this approximation
as well as the derivation of these equations to \cite[Appendix~C]{SSZZ2d}.

\subsection{Iterative PDE arguments}

We then use PDE estimate to 
obtain the following lemma. This combined together
with the above identity \eqref{e:alg-id} then yields the formula in Theorem~\ref{th:1.3}.

\begin{lemma} 
 $ \widehat{Q_N}/N\to 0$ as $N\to \infty$.
 \end{lemma}
The lemma is of course not obvious at all, since there is a sum in  \eqref{e:QN}.
One needs to gain some factors of $1/N$.
 
To explain the proof of this lemma, we only focus on the first term on the RHS
of \eqref{e:QN}. Again, we decompose $\Phi_i=Z_i+Y_i$. We have
$$
 \E \Big(\Wick{Z_i \sum_j Z_j^2 (y)}Z_i(z)\Big) =0
 $$
so it remains to bound the other terms involving $Y$.
In fact, via Holder inequality
and other manipulations, it boils down  to proving bounds of the following type  
$$
 \E \Big(\sum_i \|Y_i\|_{L^2}^2\Big)^p  < C\;,
 \qquad
  \E \|Y_i\|_{L^p}^p  < N^{\frac{p}{2}}\;.
 $$
We roughly discuss the proof for the first bound, and the second one is similar.
Let $F = \frac1N \|\sum_{i=1}^N Y_i^2\|_{L^2}^2$ and $U=\sum_{i=1}^N \|Y_i^2\|_{L^2}^2$. Similarly as how we got \eqref{e:centerD}, we can prove that
\begin{equ}[se20]
\E [U^{q-1} F ] +(m-A)  \E [U^q] \le C + C N^{-1/2} (\E U^{q+1})^{q/(q+1)}\;.
\end{equ}
Here $A$ is a constant which plays similar role as in \eqref{e:centerD}.
	The strategy now is to first use the dissipative quantity on the LHS  to obtain $\E(U^{q}) \leq CN^{\frac{q-1}{2}}$, and then use the massive term on the LHS to iteratively decrease the power of $N$ and eventually arrive at $\E(U^{q}) \leq C$.  
	
	Indeed, first observe that $F \geq N^{-1}U^{2}$ so that $\E(U^{q-1}F) \geq N^{-1}\E(U^{q+1})$.  Hence, Young's inequality with exponents $(q+1,\frac{q+1}{q})$ leads to  $\E(U^{q}) \leq CN^{\frac{q-1}{2}}$. 

%

Discarding the dissipative term, using \eqref{se20} again we have
	\begin{align}\label{szz4}
	\E  U^q \lesssim (\E  U^{q+1})^{\frac{q}{q+1}}N^{-1/2}+1
	\lesssim N^{\frac{q^2}{2(q+1)}-\frac{1}{2}}+1.
	\end{align}
We can then iterate this procedure, and eventually bring down the power of $N$ to $0$ ! See \cite[Sec.~6]{SSZZ2d} for details.

\section{Converges of (renormalized) linear $\sigma$-model to Gaussian free field in 3D}

As we explained in Sec~\ref{sec:SQ}, the analysis in three space dimensions
is significantly harder. 
Hairer's theory of regularity structures produces the first construction of local solution of dynamical $\Phi^4_3$.
Our discussion below is then much motivated by  the
SPDE construction of $\Phi^4_3$ quantum field theory by Gubinelli--Hofmanova \cite{GH21}, based on paracontrolled distributions.

We first recall some basics of paracontrolled distributions.
 In general, as mentioned in Sec.~\ref{sec:Fixed point}
 the product $fg$ of two distributions $f\in \bC^\alpha, g\in \bC^\beta$ is well defined if and only if $\alpha+\beta>0$. In terms of Littlewood-Paley blocks, the product $fg$ of two distributions $f$ and $g$ can be formally decomposed as
$$
fg=f\prec g+f\circ g+f\succ g,
$$
with \footnote{The block $\Delta_j f$ is basically the Fourier modes of order $2^j$ of $f$}
$$f\prec g=g\succ f=\sum_{j\geq-1}\sum_{i<j-1}\Delta_if\Delta_jg, \quad f\circ g=\sum_{|i-j|\leq1}\Delta_if\Delta_jg.$$
We also denote
$	\succcurlyeq \, \eqdef \, \succ+\circ$ and
$ \preccurlyeq\,\eqdef\, \prec+\circ$.
We write, for some constant $L>0$,
\begin{align}\label{zmm02}\UU_>\eqdef\sum_{j>L}\Delta_j,\quad \UU_\leq \eqdef \sum_{j\leq L}\Delta_j \;.
\end{align}

In \cite{SZZ3d}, we considered the  linear $\sigma$-model 
on three dimensional torus, which is defined as a measure $\nu^N$
 given by the limit of the following measures on the discrete torus $\mathbb{T}^3_\eps$
$$
\frac{1}{Z}\exp\Big(
-\int_{\mathbb{T}^3_\eps} \sum_{j=1}^N|\nabla_\eps \Phi_j|^2+\Big(m-\frac{N+2}N\l a_\eps+\frac{3(N+2)}{N^2}\l^2 \tilde b_\eps\Big)\sum_{j=1}^N\Phi_j^2+\frac{\l}{2N}\Big(\sum_{j=1}^N\Phi_j^2\Big)^2\dif x
\Big)\mathcal D \Phi 
$$
as lattice spacing $\eps\to 0$.
Here
 $a_\eps$ and $\tilde{b}_\eps$ are renormalization constants given below.
\footnote{These renormalization constants are the same as in 
the one-component $\Phi^4_3$ QFT, namely  $a_\eps$ diverges at rate $\eps^{-1}$ and $\tilde{b}_\eps$ diverges logarithmically.}
We showed that for sufficiently large $m$ or small $\l$, as $N \to \infty$, its $k$-component marginals converge to products of the Gaussian free field measure. 
The basic philosophy is similar as in Section~\ref{sec:GFF2d}, 
but the analysis in 3D case is much more complicated than its counterpart in 2D. 

We start with the  renormalized stochastic dynamics. Let $\xi_{i,\varepsilon}$ be a  mollification of the space-time white noise $\xi_i$  on $\R\times \mathbb{T}^3$. 
The corresponding SPDE is given by the following 
\begin{align}\label{eq:ap}
\LL \Phi_{i,\eps}
+\frac{\l}{N}\sum_{j=1}^N\Phi_{j,\eps}^2\Phi_{i,\eps}
+\Big(-\frac{N+2}{N}\l a_\eps+\frac{3(N+2)}{N^2}\l^2 \tilde{b}_\eps\Big)\Phi_{i,\eps}=\xi_{i,\eps}\;.
\end{align}
Let  $Z_i$ and $Z_{i,\varepsilon}$ be i.i.d. mean zero Gaussian processes which are
stationary solutions to
\begin{equation}\label{eq:li1}
\LL Z_i=\xi_i \;,
\qquad
\LL {Z}_{i,\varepsilon}=\xi_{i,\varepsilon}\;.
\end{equation}
We then proceed in a kind of similar manner as in the previous sections, namely, write $\Phi$ as a sum of 
(explicit, stochastic) leading order terms plus a remainder,
establish an energy identity for this remainder,
bound various terms in this identity,
and finally apply these bounds to prove the Gaussian free field limit.

\subsection{Stochastic leading terms}
\label{sec:Stochastic bounds}

As in the study of $\Phi^4_3$, we 
first introduce the ``leading order perturbative terms'',
but now these terms come with ``indices on noises''.
First, we have the Wick powers:
\begin{equ}[e:wicks]
\cZ_{ij}^{\<2>}=
\begin{cases}
\lim\limits_{\varepsilon\to0}(Z_{i,\varepsilon}^2-a_\varepsilon)  &  (i=j)\\
 \lim\limits_{\varepsilon\to0}Z_{i,\varepsilon}Z_{j,\varepsilon} & (i\neq j)
\end{cases}
\quad
\cZ_{ijj,\eps}^{\<3>} =
\begin{cases}
Z_{i,\varepsilon}^3-3a_\varepsilon Z_{i,\varepsilon}   & (i=j)\\
Z_{i,\varepsilon}Z_{j,\varepsilon}^2-a_\varepsilon Z_{i,\varepsilon} & (i\neq j)
\end{cases}
\end{equ}
 where $a_\varepsilon=\mathbf{E}[ Z_{i,\varepsilon}^2(0,0)]$.
Note that these are basically the same as in \eqref{e:wick-tilde}
except that the second definition here still depends on $\eps$ (i.e. no ``lim'')
and in 3D the constant diverges as $a_\eps \sim 1/\eps$.
Here we could think of the subscripts as indices 
assigned to the noises in the trees in the ``handwriting'' order.
The limit for $\cZ_{ij}^{\<2>}$ is in $C_T\bC^{-1-\kappa}$ for $\kappa>0$ and $T>0$.

Let $\tilde{\cZ}_{ijj,\eps}^{\<30>}$ be the stationary solution to $\LL \tilde{\cZ}_{ijj,\eps}^{\<30>}=\cZ_{ijj,\eps}^{\<3>}$, i.e. 
$$\tilde{\cZ}_{ijj,\eps}^{\<30>}=\int_{-\infty}^tP_{t-s}\cZ_{ijj,\eps}^{\<3>}(s)\dif s:=\tilde{\cI}\cZ_{ijj,\eps}^{\<3>}.$$
Set
$\tilde{\cZ}_{ijj}^{\<30>}=\lim\limits_{\varepsilon\to0}\tilde{\cZ}_{ijj,\eps}^{\<30>}$
where the limit is in  $C_T\bC^{\frac12-\kappa}$. 
To define the ``next order perturbative terms'',
let $\sD \eqdef m-\Delta$ and define
\begin{equs}[2]
\tilde{\cZ}_{ij,k\ell,\eps}^{\<22>} &\eqdef
\sD^{-1}(\cZ_{ij,\eps}^{\<2>})\circ \cZ_{k\ell,\eps}^{\<2>}- c_1 \tilde{b}_\eps \;,
\qquad
& {\cZ}_{ij,k\ell,\eps}^{\<22>} \eqdef
\cI(\cZ_{ij,\eps}^{\<2>})\circ \cZ_{k\ell,\eps}^{\<2>}- c_1 b_\eps(t) \;,
\\
\tilde{\cZ}_{ijj,k,\eps}^{\<31>} &\eqdef \tilde{\cZ}_{ijj,\eps}^{\<30>}\circ Z_{k,\eps} \;,
\qquad
&\tilde{\cZ}_{ijj,ik,\eps}^{\<32>}\eqdef
\tilde{\cZ}_{ijj,\eps}^{\<30>}\circ \cZ_{ik,\eps}^{\<2>}- c_2 \tilde{b}_\eps Z_{j,\eps}\;.
\end{equs}
Here $c_1$ equals $\frac12$ if $i=k\neq j=\ell$  or $i=\ell\neq j=k$,
equals $1$ if $i=k= j=\ell$, and is $0$ otherwise;
and $c_2 $ equals $1$ if $j=k\neq i$, equals $3$ if $j=k= i$, and equals $0$ otherwise.
Also $b_{\eps}(t)=\E[\cI(\cZ_{ii,\eps}^{\<2>})\circ \cZ_{ii,\eps}^{\<2>}]$ and  $\tilde{b}_{\eps}=\E[\sD^{-1}(\cZ_{ii,\eps}^{\<2>})\circ \cZ_{ii,\eps}^{\<2>}]$ are  renormalization constants
and $|b_\eps-\tilde{b}_\eps|\lesssim  t^{-\gamma}$  for any $\gamma>0$ uniformly in $\eps$.
We denote collectively
\begin{equ}[e:ZZZ]
\mathbb{Z}_\eps\eqdef(Z_{i,\eps},\cZ_{ij,\eps}^{\<2>},\tilde{\cZ}_{ijj,\eps}^{\<30>}, \tilde{\cZ}_{ijj,k,\eps}^{\<31>}, \tilde{\cZ}_{ij,k\ell,\eps}^{\<22>}, {\cZ}_{ij,k\ell,\eps}^{\<22>}, \tilde{\cZ}_{ijj,ik,\eps}^{\<32>}).
\end{equ}
Then, we can prove that each of these objects
converges in probability as $\eps\to 0$ in $C_T\bC^{\alpha_{\tau}}$,
where for each tree $\tau$, $\alpha_\tau = (-\frac52-\kappa)\# \mbox{noises} + 2\# \mbox{edges} $. For example for $\tau = \<2>$, $\alpha_\tau = -1-2\kappa$.

Before decomposing $\Phi$, an important idea in \cite{GH21} is to introduce 
one more auxiliary object $X_i$ defined via the following  equation:
\begin{align}\label{eq:211}
X_i=-\frac{\l}{N}\sum_{j=1}^N\Big(2\cI(X_j\prec \UU_> \cZ^{\<2>}_{ij})+\cI(X_i\prec \UU_>  \cZ^{\<2>}_{jj})+\tilde{\cZ}^{\<30>}_{ijj}\Big).
\end{align}

\subsection{Decomposition}

With the stochastic objects at hand, we have the following decompositions: 
$$\Phi_{i}=Z_i+X_i+Y_i$$ 
with $Y_i$ satisfying the following equation  
\footnote{In fact some products in \eqref{eq:YY} are understood via renormalization, see \cite{SZZ3d} for details.}
\begin{equs}\label{eq:YY}
{}\LL Y_i
&=
-\frac{\l}{N}\sum_{j=1}^N\bigg(Y_j^2Y_i+(X_j^2+2X_jY_j)(X_i+Y_i)+Y_j^2X_i+(X_j+Y_j)^2Z_i
\\
&+2(X_j+Y_j)(X_i+Y_i)Z_j+2X_j\prec \UU_\leq\cZ^{\<2>}_{ij} +X_i\prec \UU_\leq\cZ^{\<2>}_{jj}
\\
&+2Y_j\prec\cZ^{\<2>}_{ij}+Y_i\prec\cZ^{\<2>}_{jj}+2(X_j+Y_j)\succcurlyeq\cZ^{\<2>}_{ij}+(X_i+Y_i)\succcurlyeq\cZ^{\<2>}_{jj}\bigg) \;,
\\Y_i(0)=&\,\Phi_i(0)-Z_i(0)-X_i(0) \;.
\end{equs}
Here the first line in \eqref{eq:YY} are the expansion of $(Y_j+X_j)^2(Y_i+X_i+Z_i)$; the terms containing $\cZ^{\<2>}_{ij}$ and $\cZ^{\<2>}_{jj}$ correspond to the remaining terms in
 the paraproduct expansion of $2(Y_j+X_j)Z_jZ_i$ and
$(Y_i+X_i)Z_j^2$, respectively.

Recall that $\mathscr{D}=m-\Delta$ and we define
\begin{equation}\label{phi}
\varphi_i\eqdef Y_i+\mathscr{D}^{-1}\frac{\l}{N}\sum_{j=1}^N(2Y_j\prec\cZ^{\<2>}_{ij}+Y_i\prec\cZ^{\<2>}_{jj})\eqdef Y_i+\mathscr{D}^{-1}P_i,
\end{equation}
where in the last step we defined $P_i$.
We can then prove the following energy identity

\bl\label{den}
(Energy balance)
\begin{align}
\frac{1}{2} \sum_{i=1}^{N}\frac{\dif}{\dif t}\|Y_{i}\|_{L^{2}}^{2}+m\sum_{i=1}^{N}\| \varphi_{i}\|_{L^{2}}^{2}+\sum_{i=1}^{N}\|\nabla \varphi_{i}\|_{L^{2}}^{2}+\frac{\l}{N}\Big \|\sum_{i=1}^{N}Y_{i}^{2} \Big \|_{L^{2}}^{2}=\Theta+\Xi. \nonumber
\end{align}
Here
\begin{align*}
\Theta=
\sum_{i=1}^N \langle \mathscr{D}^{-1}P_i,P_i\rangle
-\frac{\l}{N}\sum_{i,j=1}^N \Big(2D(Y_i,{\cZ}^{\<2>}_{ij},Y_j)+D(Y_i,{\cZ}^{\<2>}_{jj},Y_i)\Big),
\end{align*}
and
\begin{align*}
\Xi=&-\frac{\l}{N}\sum_{i,j=1}^N\bigg\langle (X_j^2+2X_jY_j)(X_i+Y_i)+Y_j^2X_i+(X_j+Y_j)^2 Z_i
\\
&+2(X_j+Y_j)(X_i+Y_i)Z_j+2X_j\prec \UU_\leq {\cZ}^{\<2>}_{ij} +X_i\prec \UU_\leq {\cZ}^{\<2>}_{jj}
\\&+2(X_j+Y_j)\succ {{\cZ}^{\<2>}_{ij}}+2X_j\circ {{\cZ}^{\<2>}_{ij}}+(X_i+Y_i)\succ {\cZ}^{\<2>}_{jj}+X_i\circ {\cZ}^{\<2>}_{jj} \;\; ,\;\; Y_i\bigg\rangle.
\end{align*}
\el

It is straightforward to derive the above energy identity, at least formally.
In this derivation, however, 
new renormalization would appears to be necessary
when we take inner product, but in the end they cancel each other. See 
\cite{SZZ3d} for details.

The next step, similarly as in the 2D case, is to prove uniform in $N$ estimates based on  Lemma \ref{den}.

The main results are Theorems~\ref{Y:T4} and \ref{Y:T4a}.
The key step to prove these theorems
 is to obtain the following bound
\begin{equation}\label{zmm4}
\aligned
\int_0^T(\Theta+\Xi)\dif t 
&\le \delta\Big(\sum_{j=1}^N\|\nabla  \varphi_{j}\|_{L_T^2L^2}^2+\sum_{j=1}^N\| Y_{j}\|_{L_T^2H^{1-2\kappa}}^2+\frac{\l}{N}\Big\|\sum_{i=1}^NY_{i}^2\Big\|_{L^2_TL^2}^2\Big)
\\&+C_\delta\l(1+\l^{55})\int_0^T\Big(\sum_{i=1}^N\|Y_i\|_{L^2}^2\Big)(R_N^1+R_N^2+Q_N^3)\dif s+\l(1+\l) Q_N^4
\endaligned\end{equation}
for a small constant $\delta>0$
with certain quantities $R_N^1, R_N^2, Q_N^3, Q_N^4$ 
which only depend on the ``perturbative objects'' \eqref{e:ZZZ}
and $X$ in \eqref{eq:211}.
The main technical step is then to show:

\bt\label{Y:T4} There exists
a constant $C>0$ independent of  $N$, $\l$ and $m\geq1$ such that
\begin{align*}
&\Big(\sum_{j=1}^N\|Y_{j}(T)\|_{L^2}^2\Big)+\frac12\sum_{j=1}^N\|\nabla  \varphi_{j}\|_{L_T^2L^2}^2+m\sum_{j=1}^N\|Y_{j}\|_{L_T^2L^2}^2
\\
&\qquad +\frac{1}{8}\sum_{j=1}^N\| Y_{j}\|_{L_T^2H^{1-2\kappa}}^2+\frac{\l}{N}\Big\|\sum_{i=1}^NY_{i}^2\Big\|_{L^2_TL^2}^2
\\
& \leq \Big(\sum_{j=1}^N\|Y_j(0)\|_{L^2}^2\Big)+\l(1+\lambda) Q_N^4+\sum_{j=1}^N\|Y_{j}\|_{L_T^2L^2}^2
\\&\qquad+C\l(1+\l^{55})\int_0^T\Big(\sum_{i=1}^N\|Y_i\|_{L^2}^2\Big)(R_N^1+R_N^2+Q_N^3)\dif s.
\end{align*}
\et

Furthermore by using the dissipation effect from the term $\frac{1}{N}\|\sum_{i=1}^NY_{i}^2\|_{L^2_TL^2}^2$, the empirical averages of the $L^2$ norms of $Y_i$ can be controlled pathwise in terms of a $\mathbb{Z}$ dependent quantity 
 $Q(\mathbb{Z})$ with finite moment, as stated in the following theorem.

\bt\label{Y:T4a}
There exists 
 $Q(\mathbb{Z})$ with $\E Q(\mathbb{Z})\le C$  independent of $N$, such that
\begin{align*}
&\sup_{t\in[0,T]}\frac{1}{N}\sum_{j=1}^N\|Y_{j}(t)\|_{L^2}^2+\frac{m}{N}\sum_{j=1}^N\| Y_{j}\|_{L_T^2L^2}^2+\frac{1}{2N}\sum_{j=1}^N\|\nabla  \varphi_{j}\|_{L^2_TL^2}^2
\\
&\qquad\qquad\qquad\qquad
+\frac{1}{8N}\sum_{j=1}^N\| Y_{j}\|_{L_T^2H^{1-2\kappa}}^2+\frac{\l}{2N^2}\Big\|\sum_{i=1}^NY_{i}^2\Big\|_{L^2_TL^2}^2
\;
\leq Q(\mathbb{Z})+\frac{2}{N}\sum_{i=1}^N\|Y_i(0)\|_{L^2}^2.
\end{align*}
\et

The above bounds are inputs of the proof of convergence of measures
which we discuss below.

\subsection{Convergence of measures}\label{sec:GFF3d}

The proof of convergence of the renormalized linear sigma model in 3D to 
Gaussian free field is based on similar idea as in the 2D case Section~\ref{sec:GFF2d}.
However, there are several  new technical issues in 3D which we now explain.

Recall that in 3D we have a new decomposition
$$\Phi_{i}=Z_i+X_i+Y_i$$ 
where $X_i $ solves \eqref{eq:211} and 
$Y_i$ solves \eqref{eq:YY}.
Similarly as Section~\ref{sec:GFF2d} we 
consider  a  stationary process $(\Phi_i, Z_i)_{1\leq i\leq N}$ such that the components $\Phi_i$ is the limiting solution to \eqref{eq:ap} and
$Z_i$ is stationary  solution to $\LL Z_i=\xi_i$. 
Similarly as in \eqref{se8}, the convergence of measures
eventually boils down to proving
\begin{equation}\label{bdd:1}
\aligned
\E\|\tilde \Phi_i&(0)  -Z_i(0)\|_{H^{\frac{1}{2}-2\kappa}}^2
=\frac{1}{t}\int_0^t\E\|\tilde \Phi_i(s)-Z_i(s)\|_{H^{\frac{1}{2}-2\kappa}}^2 \dif s
\\
&\leq \frac{2}{t}\int_0^t\E\|X_i(s)\|_{H^{\frac{1}{2}-2\kappa}}^2 \dif s+\frac{2}{t}\int_0^t\E\|Y_i(s)\|_{H^{\frac{1}{2}-2\kappa}}^2 \dif s \stackrel{?}{\leq} CN^{-1}.
\endaligned
\end{equation}
For the  term with $X$,
it turns out that 
we have 
\begin{equ}[e:sumX]
\sum_{i=1}^N\|X_i\|_{L_T^2H^{\frac12-2\kappa}}^2  \lesssim \l^2 Q_N^0,
\end{equ}
where $Q_N^0$ only involves the ``perturbative objects''
$$
Q^0_N\eqdef\frac{1}{N^2}\sum_{i=1}^N \Big\|\sum_{j=1}^N\tilde{\cZ}^{\<30>}_{ijj}\Big\|^2_{L_T^2H^{\frac12-2\kappa}}
$$
and we have
 $\E| Q_N^0|^q\lesssim1$
 for every $q\geq 1$ uniformly in $N$.  The nontriviality of \eqref{e:sumX} is that the left hand side appears to be $O(N)$ but the right hand side is of order $O(1)$,
 despite of the fact that $Q^0_N$ appears to be ``three sums times $1/N^2\,$''.
 Thus we essentially gain a factor of $1/N$.
 \eqref{e:sumX} immediately gives 
\begin{align}\label{bdd:2}
\frac{1}{t}\int_0^t\E\|X_i(s)\|_{H^{\frac{1}{2}-2\kappa}}^2 \dif s=\frac{1}{tN}\int_0^t\E\sum_{i=1}^N\|X_i(s)\|_{H^{\frac{1}{2}-2\kappa}}^2 \dif s
\lesssim \frac{1}{tN}\E Q_N^0\lesssim \frac{C_t}{tN},
\end{align}
so we get the desired estimate for the first term in \eqref{bdd:1}.

 The proof of  \eqref{e:sumX}  is based on a priori estimates 
 on the PDE  \eqref{eq:211}. We only 
 demonstrate why  the expectation of $Q_N^0$ is bounded uniformly in $N$.
  We have for $s=\frac12-2\kappa$
	\begin{align*}
	\E\frac{1}{N^2} \sum_{i=1}^{N}  \Big\|\sum_{j=1}^{N} \nabla^{s}\tilde{\cZ}^{\<30>}_{ijj}\Big\|^2_{L_T^2L^2}
	=
	\frac{1}{N^2}\sum_{i,j_1,j_2=1}^N\E \Big\la\nabla^{s} \tilde{\cZ}^{\<30>}_{ij_1j_1},\nabla^{s} \tilde{\cZ}^{\<30>}_{ij_2j_2}\Big\ra_{L^2_TL^2} \sim O(1).
	\end{align*}
	We have 3 summation indices and a factor $1/N^2$.
	 However, if $i,j_1,j_2$ are all different, by independence
	and the fact that Wick products are mean zero, the terms are zero. (Note that 
this comes from the same idea as \eqref{eq:R0} in 2D case.)

We remark
that the 
same mechanism also applies to higher order perturbative objects such as
$\tilde{\cZ}^{\<32>}$ (see \cite[Lemma~3.4]{SZZ3d}).

Next, we turn to the term with $Y$ in \eqref{bdd:1}. The start-point is similar 
with \eqref{e:centerD}, but the analysis is more involved.
From Theorem~\ref{Y:T4}
we can deduce a bound of the form
\begin{align}
&\frac18\int_0^t\E\Big(\sum_{i=1}^N\|Y_i(s)\|_{H^{1-2\kappa}}^2\Big) \dif s+(m-1)\int_0^t\E\Big(\sum_{i=1}^N\|Y_i(s)\|_{L^2}^2\Big) \dif s+\frac{\l}{N}\E\int_0^t\Big\|\sum_{i=1}^NY_i^2\Big\|_{L^2}^2\dif s\no
\\
&\leq 2N\E\|\tilde{\Phi}_i(0)-Z_i(0)\|_{L^2}^2+2\sum_{i=1}^N\E\|X_i(0)\|_{L^2}^2
+C_t
+C \l 
\int_0^t\E\Big(\sum_{i=1}^N\|Y_i(s)\|_{L^2}^2\Big)\E R_N \dif s  \label{sz1}
\\
&\quad
+C\l \int_0^t\E\Big[\Big(\sum_{i=1}^N\|Y_i(s)\|_{L^2}^2\Big)
\Big| R_N -\E[R_N]\Big|\Big] \dif s. \no
\end{align}
Here $C$ is independent of $\l, m$ and $N$.

The quantity $R_N$ has a complicated expression, which we do not give here,
but the only point which matters is that it is defined via the perturbative objects
in Section~\ref{sec:Stochastic bounds}.
Moreover,
all the terms in
 $R_N - \E R_N$
are  summations of terms of the form
$$
 \frac{1}{N^l} \sum_{i_1\cdots i_l=1}^N M_{i_1,\cdots,i_l}
$$
for certain choice of $l$,
where each $M_{i_1,\cdots,i_l}$ is mean-zero, has bounded second moment,
 and
they satisfy a crucial independence condition:
$M_{i_1,\cdots,i_l} $ and $M_{j_1,\cdots,j_l}$ are independent when
the $2l$ indices $i_1,\cdots,i_l,j_1,\cdots,j_l$ are all different. 
For any such sum, one can prove that 
$$
\E \Big[ \Big( \frac{1}{N^l} \sum_{i_1\cdots i_l=1}^N M_{i_1,\cdots,i_l} \Big)^2\Big] \le C/N
$$
where $C$ only depends on $l$ and is independent of $N$.
From these estimates,
we can eventually deduce \eqref{bdd:1}, see \cite[Sec.~5]{SZZ3d} for details.

We remark that we can also show tightness of some observables
in suitable Besov spaces in 3D, see \cite[Theorem~1.2]{SZZ3d}.

\section{Discussions on possible future directions}

We complete this review by discussing some possible future problems.
\begin{enumerate}
\item
For the linear $\sigma$-model on three dimensional torus,
we only studied stationary solutions and their large $N$  limit in \cite{SZZ3d}.
It would be interesting to generalize the mean field limit for the dynamics in 2D discussed in  Section~\ref{sec:MF2D} to 3D.
\item
In 2D, besides tightness of the two observables 
in \eqref{e:twoObs} 
 we also prove exact formulas for their correlations as $N\to \infty$ (Theorem~\ref{th:1.3});
it would be interesting to prove exact correlation formulas for other $O(N)$-invariant observables in the large $N$ limit.
In 3D, we only obtained  tightness of observables $\frac{1}{\sqrt{N}} \sum_{i=1}^N\Wick{\Phi^2_i}$ as random fields when $N$ becomes large in \cite{SZZ3d}, so
it would be interesting to derive the exact correlation formulas for it in the large $N$ limit.

\item
It would be interesting to obtain a $1/N$ expansion for correlations of the field
or of the observables, and prove bounds on the remainder of this $1/N$ expansion. In 2D, one might be able to incorporate 
the dual field approach 
discussed in Section~\ref{sec:DS} as in Kupiainen \cite[Section~4]{MR578040},
and then, bound the remainder of this $1/N$ expansion
using PDE estimates.  In fact, for the usual perturbation expansion for $\Phi^4_2$, we have already given new bounds on the remainder (which implies asymptoticity of the expansion) in \cite{SZZperturbation}.
\item
Could any results discussed in this article be generalized to infinite volume?
\item
For Fermionic models, could we make any prediction by physicists  e.g.  \cite{gross1974dynamical} rigorous? Recently, \cite{grassmannian2020} proposed a rigorous framework for stochastic quantization of fermionic fields (at least on lattice), which might be a reasonable start-point.
\item
It would be interesting to study large $N$ behavior of matrix-valued dynamics,
or develop estimates on the dynamics to deduce  large $N$ results on matrix-valued QFT (e.g. \cite{gopakumar1995mastering}), for instance \eqref{e:Phi_i-measure} with $\Phi$ being matrix valued, either on lattice or in continuum.  
Besides vector and matrix models, tensor models and their large $N$ problems were also proposed in physics.
\item
Among the matrix models, Yang--Mills models are of great interest, see the recent progress on large N problems for lattice  Yang--Mills  in 
 \cite{MR3919447}, \cite{chatterjee2016}. In equilibrium, Wilson loop observables 
 converge to deterministic limit i.e. master fields.
It would be interesting to use the dynamic to study such questions:
in \cite{stochastic-lattice-YM-2022} deterministic large $N$ limit was shown using the dynamic, and it would be interesting to study higher order $1/N$ corrections.
Another possible questions is to obtain certain deterministic limit for the dynamics as $N\to \infty$. 
The paper \cite{SSZloop} is also an initial attempt where the master loop equation for finite $N$ is derived using the stochastic dynamic on lattice. 
On two and three dimensional continuum torus, these dynamics are constructed for short time \cite{CCHS20,CCHS3D}.
\end{enumerate}

\subsection*{Acknowledgments}
 H.S. gratefully acknowledges financial support from NSF grants DMS-1954091 and CAREER DMS-2044415, and collaborations with Scott Smith, Rongchan Zhu and Xiangchan Zhu.

\subsection*{Data availability}
Data sharing is not applicable to this article as no new data were created or analyzed in this study.

\bibliographystyle{plain}
\bibliography{Reference}

\begin{thebibliography}{10}

\bibitem{abbott1976bound}
LF~Abbott, JS~Kang, and HJ~Schnitzer.
\newblock Bound states, tachyons, and restoration of symmetry in the {$1/N$}
  expansion.
\newblock {\em Physical Review D}, 13(8):2212, 1976.

\bibitem{MR1113223}
S.~Albeverio and M.~R\"{o}ckner.
\newblock Stochastic differential equations in infinite dimensions: solutions
  via {D}irichlet forms.
\newblock {\em Probab. Theory Related Fields}, 89(3):347--386, 1991.

\bibitem{grassmannian2020}
Sergio Albeverio, Luigi Borasi, Francesco~C De~Vecchi, and Massimiliano
  Gubinelli.
\newblock Grassmannian stochastic analysis and the stochastic quantization of
  euclidean fermions.
\newblock {\em arXiv preprint arXiv:2004.09637}, 2020.

\bibitem{AK17}
Sergio Albeverio and Seiichiro Kusuoka.
\newblock The invariant measure and the flow associated to the $\phi^4_3
  $-quantum field model.
\newblock {\em Annali della Scuola Normale Superiore di Pisa-Classe di
  Scienze}, To appear.

\bibitem{Alfaro1983}
J.~Alfaro.
\newblock Stochastic quantization and the large-{N} reduction of {U(N)} gauge
  theory.
\newblock {\em Physical Review D}, 28(4):1001, 1983.

\bibitem{AlfaroSakita}
J.~Alfaro and B.~Sakita.
\newblock Derivation of quenched momentum prescription by means of stochastic
  quantization.
\newblock {\em Physics Letters B}, 121(5):339--344, 1983.

\bibitem{MR2864481}
Michael Anshelevich and Ambar~N. Sengupta.
\newblock Quantum free {Y}ang-{M}ills on the plane.
\newblock {\em J. Geom. Phys.}, 62(2):330--343, 2012.

\bibitem{MR2768550}
Hajer Bahouri, Jean-Yves Chemin, and Rapha\"{e}l Danchin.
\newblock {\em Fourier analysis and nonlinear partial differential equations},
  volume 343 of {\em Grundlehren der Mathematischen Wissenschaften [Fundamental
  Principles of Mathematical Sciences]}.
\newblock Springer, Heidelberg, 2011.

\bibitem{BCD18}
Isma\"{e}l Bailleul, R\'{e}mi Catellier, and Francois Delarue.
\newblock Solving mean field rough differential equations.
\newblock {\em Electron. J. Probab.}, 25:Paper No. 21, 51, 2020.

\bibitem{bardeen1976phase}
William~A Bardeen, Benjamin~W Lee, and Robert~E Shrock.
\newblock Phase transition in the nonlinear $\sigma$ model in a (2+
  $\varepsilon$)-dimensional continuum.
\newblock {\em Physical Review D}, 14(4):985, 1976.

\bibitem{baxter2016exactly}
Rodney~J Baxter.
\newblock {\em Exactly solved models in statistical mechanics}.
\newblock Elsevier, 2016.

\bibitem{berlin1952spherical}
Theodore~H Berlin and Mark Kac.
\newblock The spherical model of a ferromagnet.
\newblock {\em Physical Review}, 86(6):821, 1952.

\bibitem{MR661137}
C.~Billionnet and P.~Renouard.
\newblock Analytic interpolation and {B}orel summability of the {$({\lambda
  \over N}\Phi _{N}^{:4})_{2}$} models. {I}. {F}inite volume approximation.
\newblock {\em Comm. Math. Phys.}, 84(2):257--295, 1982.

\bibitem{brezin1973critical}
E~Br{\'e}zin and DJ~Wallace.
\newblock Critical behavior of a classical {H}eisenberg ferromagnet with many
  degrees of freedom.
\newblock {\em Physical Review B}, 7(5):1967, 1973.

\bibitem{LargeN1993}
Edouard Brezin and Spenta~R Wadia.
\newblock {\em The large {N} expansion in quantum field theory and statistical
  physics: from spin systems to 2-dimensional gravity}.
\newblock World scientific, 1993.

\bibitem{MR496278}
David Brydges and Paul Federbush.
\newblock A lower bound for the mass of a random {G}aussian lattice.
\newblock {\em Comm. Math. Phys.}, 62(1):79--82, 1978.

\bibitem{MR3299600}
Thomas Cass and Terry Lyons.
\newblock Evolving communities with individual preferences.
\newblock {\em Proc. Lond. Math. Soc. (3)}, 110(1):83--107, 2015.

\bibitem{MR3846835}
R\'{e}mi Catellier and Khalil Chouk.
\newblock Paracontrolled distributions and the 3-dimensional stochastic
  quantization equation.
\newblock {\em Ann. Probab.}, 46(5):2621--2679, 2018.

\bibitem{CCHS20}
Ajay Chandra, Ilya Chevyrev, Martin Hairer, and Hao Shen.
\newblock Langevin dynamic for the 2{D} {Y}ang-{M}ills measure.
\newblock {\em arXiv:2006.04987, To appear in Publ.Math.IH\'ES}, 2020.

\bibitem{CCHS3D}
Ajay Chandra, Ilya Chevyrev, Martin Hairer, and Hao Shen.
\newblock Stochastic quantisation of {Yang-Mills-Higgs in 3D}.
\newblock {\em arXiv preprint arXiv:2201.03487}, 2022.

\bibitem{MR3919447}
Sourav Chatterjee.
\newblock Rigorous solution of strongly coupled {$SO(N)$} lattice gauge theory
  in the large {$N$} limit.
\newblock {\em Comm. Math. Phys.}, 366(1):203--268, 2019.

\bibitem{chatterjee2016}
Sourav Chatterjee and Jafar Jafarov.
\newblock The {$1/N $} expansion for {SO(N)} lattice gauge theory at strong
  coupling.
\newblock {\em arXiv preprint arXiv:1604.04777}, 2016.

\bibitem{CDFM18}
Michele Coghi, Jean-Dominique Deuschel, Peter~K Friz, and Mario Maurelli.
\newblock Pathwise {McKean--Vlasov} theory with additive noise.
\newblock {\em Annals of Applied Probability}, 30(5):2355--2392, 2020.

\bibitem{coleman1988}
Sidney Coleman.
\newblock {\em Aspects of symmetry: selected {E}rice lectures}.
\newblock Cambridge University Press, 1988.

\bibitem{coleman1974}
Sidney Coleman, Roman Jackiw, and Hugh~David Politzer.
\newblock Spontaneous symmetry breaking in the {O(N)} model for large {N}.
\newblock {\em Physical Review D}, 10(8):2491, 1974.

\bibitem{DD03}
Giuseppe Da~Prato and Arnaud Debussche.
\newblock Strong solutions to the stochastic quantization equations.
\newblock {\em Ann. Probab.}, 31(4):1900--1916, 2003.

\bibitem{d1979confinement}
A~d'Adda, P~Di~Vecchia, and Martin L{\"u}scher.
\newblock Confinement and chiral symmetry breaking in {$CP^{n- 1}$} models with
  quarks.
\newblock {\em Nuclear Physics B}, 152(1):125--144, 1979.

\bibitem{d19781n}
A~d'Adda, M~L{\"u}scher, and P~Di~Vecchia.
\newblock A $1/n$ expandable series of non-linear $\sigma$ models with
  instantons.
\newblock {\em Nuclear Physics B}, 146(1):63--76, 1978.

\bibitem{MR3554890}
Antoine Dahlqvist.
\newblock Free energies and fluctuations for the unitary {B}rownian motion.
\newblock {\em Comm. Math. Phys.}, 348(2):395--444, 2016.

\bibitem{damgaard1987}
Poul~H Damgaard and Helmuth H{\"u}ffel.
\newblock Stochastic quantization.
\newblock {\em Physics Reports}, 152(5-6):227--398, 1987.

\bibitem{MR3982691}
Bruce~K. Driver.
\newblock A functional integral approaches to the {M}akeenko-{M}igdal
  equations.
\newblock {\em Comm. Math. Phys.}, 370(1):49--116, 2019.

\bibitem{MR3631396}
Bruce~K. Driver, Franck Gabriel, Brian~C. Hall, and Todd Kemp.
\newblock The {M}akeenko-{M}igdal equation for {Y}ang-{M}ills theory on compact
  surfaces.
\newblock {\em Comm. Math. Phys.}, 352(3):967--978, 2017.

\bibitem{MR3613519}
Bruce~K. Driver, Brian~C. Hall, and Todd Kemp.
\newblock Three proofs of the {M}akeenko-{M}igdal equation for {Y}ang-{M}ills
  theory on the plane.
\newblock {\em Comm. Math. Phys.}, 351(2):741--774, 2017.

\bibitem{MR3160067}
Weinan E and Hao Shen.
\newblock Mean field limit of a dynamical model of polymer systems.
\newblock {\em Sci. China Math.}, 56(12):2591--2598, 2013.

\bibitem{MR2680421}
L\'{a}szl\'{o} Erd\H{o}s, Benjamin Schlein, and Horng-Tzer Yau.
\newblock Derivation of the {G}ross-{P}itaevskii equation for the dynamics of
  {B}ose-{E}instein condensate.
\newblock {\em Ann. of Math. (2)}, 172(1):291--370, 2010.

\bibitem{MR678004}
J.~Fr\"{o}hlich, A.~Mardin, and V.~Rivasseau.
\newblock Borel summability of the {$1/N$} expansion for the {$N$}-vector
  [{${\rm O}(N)$} nonlinear {$\sigma $}] models.
\newblock {\em Comm. Math. Phys.}, 86(1):87--110, 1982.

\bibitem{MR436829}
J.~Fr\"{o}hlich and B.~Simon.
\newblock Pure states for general {$P(\phi )_{2}$} theories: construction,
  regularity and variational equality.
\newblock {\em Ann. of Math. (2)}, 105(3):493--526, 1977.

\bibitem{MR887102}
James Glimm and Arthur Jaffe.
\newblock {\em Quantum physics}.
\newblock Springer-Verlag, New York, second edition, 1987.
\newblock A functional integral point of view.

\bibitem{MR3468297}
Francois Golse.
\newblock On the dynamics of large particle systems in the mean field limit.
\newblock In {\em Macroscopic and large scale phenomena: coarse graining, mean
  field limits and ergodicity}, volume~3 of {\em Lect. Notes Appl. Math.
  Mech.}, pages 1--144. Springer, [Cham], 2016.

\bibitem{gopakumar1995mastering}
Rajesh Gopakumar and David~J Gross.
\newblock Mastering the master field.
\newblock {\em Nuclear Physics B}, 451(1-2):379--415, 1995.

\bibitem{gross1974dynamical}
David~J Gross and Andre Neveu.
\newblock Dynamical symmetry breaking in asymptotically free field theories.
\newblock {\em Physical Review D}, 10(10):3235, 1974.

\bibitem{GHglobal}
Massimiliano Gubinelli and Martina Hofmanov\'{a}.
\newblock Global solutions to elliptic and parabolic {$\Phi^4$} models in
  {E}uclidean space.
\newblock {\em Comm. Math. Phys.}, 368(3):1201--1266, 2019.

\bibitem{GH21}
Massimiliano Gubinelli and Martina Hofmanov\'{a}.
\newblock A {PDE} construction of the {E}uclidean {$\phi_3^4$} quantum field
  theory.
\newblock {\em Comm. Math. Phys.}, 384(1):1--75, 2021.

\bibitem{guionnet2019asymptotics}
Alice Guionnet.
\newblock {\em Asymptotics of random matrices and related models: the uses of
  {Dyson-Schwinger} equations}, volume 130.
\newblock American Mathematical Soc., 2019.

\bibitem{Hairer14}
M.~Hairer.
\newblock A theory of regularity structures.
\newblock {\em Invent. Math.}, 198(2):269--504, 2014.

\bibitem{MR1686539}
K.~R. Ito and H.~Tamura.
\newblock {$N$} dependence of upper bounds of critical temperatures of {${\rm
  2D\ O}(N)$} spin models.
\newblock {\em Comm. Math. Phys.}, 202(1):127--168, 1999.

\bibitem{MR3317577}
Pierre-Emmanuel Jabin.
\newblock A review of the mean field limits for {V}lasov equations.
\newblock {\em Kinet. Relat. Models}, 7(4):661--711, 2014.

\bibitem{Jafar}
Jafar Jafarov.
\newblock Wilson loop expectations in {$ SU (N) $} lattice gauge theory.
\newblock {\em arXiv preprint arXiv:1610.03821}, 2016.

\bibitem{MR1773042}
Arthur Jaffe.
\newblock Constructive quantum field theory.
\newblock In {\em Mathematical physics 2000}, pages 111--127. Imp. Coll. Press,
  London, 2000.

\bibitem{kac1971spherical}
Mark Kac and Colin~J Thompson.
\newblock Spherical model and the infinite spin dimensionality limit.
\newblock {\em Physica Norvegica}, 5(3-4):163--168, 1971.

\bibitem{MR1465436}
Gopinath Kallianpur and Jie Xiong.
\newblock {\em Stochastic differential equations in infinite-dimensional
  spaces}, volume~26 of {\em Institute of Mathematical Statistics Lecture
  Notes---Monograph Series}.
\newblock Institute of Mathematical Statistics, Hayward, CA, 1995.

\bibitem{kazakov1981wilson}
Vladimir~A Kazakov.
\newblock Wilson loop average for an arbitrary contour in two-dimensional $u
  (n)$ gauge theory.
\newblock {\em Nuclear Physics B}, 179(2):283--292, 1981.

\bibitem{kazakov1980non}
Vladimir~A Kazakov and Ivan~K Kostov.
\newblock Non-linear strings in two-dimensional {$U(\infty)$} gauge theory.
\newblock {\em Nuclear Physics B}, 176(1):199--215, 1980.

\bibitem{MR1686543}
C.~Kopper.
\newblock Mass generation in the large {$N$}-nonlinear {$\sigma$}-model.
\newblock {\em Comm. Math. Phys.}, 202(1):89--126, 1999.

\bibitem{MR1328264}
C.~Kopper, J.~Magnen, and V.~Rivasseau.
\newblock Mass generation in the large {$N$} {G}ross-{N}eveu-model.
\newblock {\em Comm. Math. Phys.}, 169(1):121--180, 1995.

\bibitem{MR582622}
A.~Kupiainen.
\newblock {$1/n$} expansion---some rigorous results.
\newblock In {\em Mathematical problems in theoretical physics ({P}roc.
  {I}nternat. {C}onf. {M}ath. {P}hys., {L}ausanne, 1979)}, volume 116 of {\em
  Lecture Notes in Phys.}, pages 208--210. Springer, Berlin-New York, 1980.

\bibitem{MR578040}
Antti~J. Kupiainen.
\newblock {$1/n$} expansion for a quantum field model.
\newblock {\em Comm. Math. Phys.}, 74(3):199--222, 1980.

\bibitem{MR574175}
Antti~J. Kupiainen.
\newblock On the {$1/n$} expansion.
\newblock {\em Comm. Math. Phys.}, 73(3):273--294, 1980.

\bibitem{MR2295621}
Jean-Michel Lasry and Pierre-Louis Lions.
\newblock Mean field games.
\newblock {\em Jpn. J. Math.}, 2(1):229--260, 2007.

\bibitem{MR2006374}
Thierry L\'{e}vy.
\newblock Yang-{M}ills measure on compact surfaces.
\newblock {\em Mem. Amer. Math. Soc.}, 166(790):xiv+122, 2003.

\bibitem{MR2667871}
Thierry L\'{e}vy.
\newblock Two-dimensional {M}arkovian holonomy fields.
\newblock {\em Ast\'{e}risque}, (329):172, 2010.

\bibitem{Levy11}
Thierry L\'{e}vy.
\newblock The master field on the plane.
\newblock {\em Ast\'{e}risque}, (388):ix+201, 2017.

\bibitem{makeenko1979exact}
Yu~M Makeenko and Alexander~A Migdal.
\newblock Exact equation for the loop average in multicolor {QCD}.
\newblock {\em Physics Letters B}, 88(1-2):135--137, 1979.

\bibitem{MR0233437}
H.~P. McKean.
\newblock Propagation of chaos for a class of non-linear parabolic equations.
\newblock In {\em Stochastic {D}ifferential {E}quations ({L}ecture {S}eries in
  {D}ifferential {E}quations, {S}ession 7, {C}atholic {U}niv., 1967)}, pages
  41--57. Air Force Office Sci. Res., Arlington, Va., 1967.

\bibitem{moinat2020space}
Augustin Moinat and Hendrik Weber.
\newblock Space-time localisation for the dynamic $\phi^4_3$ model.
\newblock {\em Communications on Pure and Applied Mathematics},
  73(12):2519--2555, 2020.

\bibitem{moshe2003}
Moshe Moshe and Jean Zinn-Justin.
\newblock Quantum field theory in the large {N} limit: A review.
\newblock {\em Physics Reports}, 385(3-6):69--228, 2003.

\bibitem{MW18}
Jean-Christophe Mourrat and Hendrik Weber.
\newblock The dynamic {$\Phi^4_3$} model comes down from infinity.
\newblock {\em Comm. Math. Phys.}, 356(3):673--753, 2017.

\bibitem{MW17}
Jean-Christophe Mourrat and Hendrik Weber.
\newblock Global well-posedness of the dynamic {$\Phi^4$} model in the plane.
\newblock {\em Ann. Probab.}, 45(4):2398--2476, 2017.

\bibitem{ParisiWu}
Giorgio Parisi and Yong~Shi Wu.
\newblock Perturbation theory without gauge fixing.
\newblock {\em Sci. Sinica}, 24(4):483--496, 1981.

\bibitem{MR471851}
Paul~A. Pearce and Colin~J. Thompson.
\newblock The spherical limit for {$n$}-vector correlations.
\newblock {\em J. Statist. Phys.}, 17(4):189--196, 1977.

\bibitem{MR1402248}
Michael~E. Peskin and Daniel~V. Schroeder.
\newblock {\em An introduction to quantum field theory}.
\newblock Addison-Wesley Publishing Company, Advanced Book Program, Reading,
  MA, 1995.
\newblock Edited and with a foreword by David Pines.

\bibitem{polyakov1980gauge}
Alexander~M Polyakov.
\newblock Gauge fields as rings of glue.
\newblock {\em Nuclear Physics B}, 164:171--188, 1980.

\bibitem{MR1346931}
Ambar Sengupta.
\newblock Gauge theory on compact surfaces.
\newblock {\em Mem. Amer. Math. Soc.}, 126(600):viii+85, 1997.

\bibitem{MR2757706}
Ambar~N. Sengupta.
\newblock The large-{$N$} {Y}ang-{M}ills field on the plane and free noise.
\newblock In {\em Geometric methods in physics}, volume 1079 of {\em AIP Conf.
  Proc.}, pages 121--132. Amer. Inst. Phys., Melville, NY, 2008.

\bibitem{MR982418}
M.~V. Shcherbina.
\newblock The spherical limit of {$n$}-vector correlations.
\newblock {\em Teoret. Mat. Fiz.}, 77(3):460--471, 1988.

\bibitem{SSZloop}
Hao Shen, Scott~A Smith, and Rongchan Zhu.
\newblock A new derivation of the finite $ n $ master loop equation for lattice
  {Yang--Mills}.
\newblock {\em arXiv preprint arXiv:2202.00880}, 2022.

\bibitem{SSZZ2d}
Hao Shen, Scott~A Smith, Rongchan Zhu, and Xiangchan Zhu.
\newblock {Large $N$ limit of the $O (N)$ linear sigma model via stochastic
  quantization}.
\newblock {\em The Annals of Probability}, 50(1):131--202, 2022.

\bibitem{SZZperturbation}
Hao Shen, Rongchan Zhu, and Xiangchan Zhu.
\newblock {An SPDE approach to perturbation theory of $\Phi^4_2 $:
  asymptoticity and short distance behavior}.
\newblock {\em arXiv preprint arXiv:2108.11312}, 2021.

\bibitem{SZZ3d}
Hao Shen, Rongchan Zhu, and Xiangchan Zhu.
\newblock {Large $ N $ limit of the $ O (N) $ linear sigma model in 3D}.
\newblock {\em arXiv:2102.02628, To appear in Comm. Math. Phys.}, 2021.

\bibitem{stochastic-lattice-YM-2022}
Hao Shen, Rongchan Zhu, and Xiangchan Zhu.
\newblock A stochastic analysis approach to lattice {Yang--Mills} at strong
  coupling.
\newblock {\em arXiv preprint arXiv:2204.12737}, 2022.

\bibitem{MR1373007}
I.~M. Singer.
\newblock On the master field in two dimensions.
\newblock In {\em Functional analysis on the eve of the 21st century, {V}ol. 1
  ({N}ew {B}runswick, {NJ}, 1993)}, volume 131 of {\em Progr. Math.}, pages
  263--281. Birkh\"{a}user Boston, Boston, MA, 1995.

\bibitem{spohn2012large}
H.~Spohn.
\newblock {\em Large Scale Dynamics of Interacting Particles}.
\newblock Texts and Monographs in Physics. Springer, 1991.

\bibitem{stanley1968spherical}
Harry~Eugene Stanley.
\newblock Spherical model as the limit of infinite spin dimensionality.
\newblock {\em Physical Review}, 176(2):718, 1968.

\bibitem{symanzik1977}
K~Symanzik.
\newblock 1/{N} expansion in {$P(\varphi^2)_{4-\epsilon}$} theory {I}. massless
  theory $0<\epsilon<2$.
\newblock {\em preprint {DESY}}, 77(05), 1977.

\bibitem{MR1108185}
Alain-Sol Sznitman.
\newblock Topics in propagation of chaos.
\newblock In {\em \'{E}cole d'\'{E}t\'{e} de {P}robabilit\'{e}s de
  {S}aint-{F}lour {XIX}---1989}, volume 1464 of {\em Lecture Notes in Math.},
  pages 165--251. Springer, Berlin, 1991.

\bibitem{tHooft1974planar}
G~t'Hooft.
\newblock A planar diagram theory for strong interactions.
\newblock {\em Nuclear Physics. B}, 72(3):461--473, 1974.

\bibitem{wilson1973quantum}
Kenneth~G Wilson.
\newblock Quantum field-theory models in less than 4 dimensions.
\newblock {\em Physical Review D}, 7(10):2911, 1973.

\bibitem{witten1980}
Edward Witten.
\newblock The 1/{N} expansion in atomic and particle physics.
\newblock In {\em Recent developments in gauge theories}, pages 403--419.
  Springer, 1980.

\end{thebibliography}

\end{document}